\documentclass[preprint,authoryear, fleqn, 12pt]{article}
\usepackage{amssymb}
\usepackage{mathrsfs}
\usepackage{ifpdf}
\usepackage{graphicx}
\ifpdf  \fi
\usepackage{pstricks}
\usepackage{subfigure}
\usepackage{rotating}
\usepackage{amsfonts,amssymb,slashbox}
\usepackage{amsmath}
\usepackage{mathptmx,helvet,courier,makeidx,multicol,footmisc}
\usepackage{natbib}
\usepackage{epstopdf}
\usepackage{multicol}
\usepackage{multirow}
\usepackage{setspace}
\usepackage{color,soul}
\usepackage{booktabs}
\usepackage[bookmarksnumbered=true, pdfauthor={Qi-jian Gan}]{hyperref}

\usepackage{array}

\usepackage{color}

\newcommand{\commentout}[1]{}

\newcommand{\ba}{\begin{array}}
        \newcommand{\ea}{\end{array}}
\newcommand{\bc}{\begin{center}}
        \newcommand{\ec}{\end{center}}
\newcommand{\bdm}{\begin{displaymath}}
        \newcommand{\edm}{\end{displaymath}}
\newcommand{\bds} {\begin{description}}
        \newcommand{\eds} {\end{description}}
\newcommand{\ben}{\begin{enumerate}}
        \newcommand{\een}{\end{enumerate}}
\newcommand{\beq}{\begin{equation}}
        \newcommand{\eeq}{\end{equation}}
\newcommand{\bfg} {\begin{figure}[h!]}
        \newcommand{\efg} {\end{figure}}
\newcommand{\bi} {\begin {itemize}}
        \newcommand{\ei} {\end {itemize}}
\newcommand{\bqn}{\begin{eqnarray*}}
        \newcommand{\eqn}{\end{eqnarray}}
\newcommand{\bqs}{\begin{eqnarray}}
        \newcommand{\eqs}{\end{eqnarray}}
\newcommand{\bsl} {\begin{slide}[8.8in,6.7in]}
        \newcommand{\esl} {\end{slide}}
\newcommand{\bss} {\begin{slide*}[9.3in,6.7in]}
        \newcommand{\ess} {\end{slide*}}
\newcommand{\btb} {\begin {table}}
        \newcommand{\etb} {\end {table}}
\newcommand{\bvb}{\begin{verbatim}}
        \newcommand{\evb}{\end{verbatim}}

\newcommand{\cas}[1]{{{\left \{ \ba #1 \ea \right. }}}

\newcommand{\refe}[1] {{(\ref {#1})}}

\def\pmb#1{\setbox0=\hbox{$#1$}%
   \kern-.025em\copy0\kern-\wd0
   \kern.05em\copy0\kern-\wd0
   \kern-.025em\raise.0433em\box0 }

\def\eop{{\hfill $\blacksquare$}}

\newtheorem{theorem}{Theorem}[section]
\newtheorem{definition}[theorem]{Definition}
\newtheorem{lemma}[theorem]{Lemma}

\begin{document}

\title{Analysis of traffic statics and dynamics in a signalized double-ring network: A Poincar\'{e} map approach}

\author{ Qi-Jian Gan\footnote{Department of Civil and Environmental Engineering, Institute of Transportation Studies, University of California, Irvine, CA 92697-3600. Email: qgan@uci.edu} and Wen-Long Jin \footnote{Department of Civil and Environmental Engineering, California Institute for Telecommunications and Information Technology, Institute of Transportation Studies, 4000 Anteater Instruction and Research Bldg, University of California, Irvine, CA 92697-3600. Email: wjin@uci.edu. Corresponding author} and Vikash V. Gayah \footnote{The Pensilvania State University, 223A Sackett Building, University Park, PA, 16801. Email: gayah@engr.psu.edu}}

\maketitle
 
\begin{abstract}
Understanding traffic statics and dynamics in urban networks is critical to develop effective control and management strategies. In this paper, we provide a novel approach to study the traffic statics and dynamics in a signalized double-ring network, which can provide insights into the operation of more general signalized traffic networks. Under the framework of the link queue model (LQM) and the assumption of a triangular traffic flow fundamental diagram, the signalized double-ring network is studied as a switched affine system. Due to periodic signal regulations, periodic density evolution orbits are formed and defined as stationary states. A Poincar\'{e} map approach is introduced to analyze the properties of such stationary states. With short cycle lengths, closed-form Poincar\'{e} maps are derived. Stationary states and their stability properties are obtained by finding and analyzing the fixed points on the Poincar\'{e} maps. It is found that a stationary state can be asymptotically stable, Lyapunov stable, or unstable. The impacts of retaining ratios and initial densities on the macroscopic fundamental diagrams (MFDs) and the gridlock times are analyzed. Multivaluedness and gridlock phenomena as well as the unstable branch with non-zero average network flow-rates are observed on the MFDs. With long cycle lengths, fixed points on the Poincar\'{e} maps are solved numerically, and the obtained stationary states and the MFDs are very similar to those with short cycle lengths.

Compared with earlier studies, this paper provides an analytical framework that can be used to provide complete and closed-form solutions to the statics and dynamics of double-ring networks. This can lead to a better understanding of how the combination of signalized intersections and turning maneuvers is expected to impact network properties, like the MFD.
\end{abstract}


\section{Introduction}\label{introduction}
The rapid increase in travel demands has led to severe congestion problems in urban transportation networks \citep{papageorgiou2003review}. This congestion has many significant negative impacts to society, including lost productivity, wasted fuel, and the creation of environmental and health problems \citep{schrank2012tti}. One recommended strategy to alleviate urban traffic congestion is improving signal timings at intersections where the current technology is deficient \citep{sorensen2008moving}. Many control strategies have been proposed since one of the earliest attempts in \citep{webster1958traffic}, but most are developed for under-saturated conditions and are not applicable to over-saturated conditions \citep{papageorgiou2003review}. In order to develop effective control and management strategies for urban networks, a fundamental challenge is to understand the static and dynamic properties of traffic in signalized networks under both under- and over-saturated conditions and a variety of other factors, e.g., demand patterns, signal settings, and route choice behaviors. 

To evaluate the impact of changes (e.g., re-timing the signals) on urban networks, one practical approach is via the macroscopic fundamental diagram (MFD)\citep{daganzo2008analytical}. The MFD is a network-level flow-density relation, which was first proposed by \citep{godfrey1969mechanism}. Recent studies have shown that MFD exists in urban networks \citep{geroliminis2008existence} and is related to many factors, such as spatial variability in congestion and density distributions \citep{buisson2009exploring, ji2010investigating, geroliminis2011properties,daganzo2011macroscopic-b}, route choices \citep{knoop2012routing}, signal settings and turning movements \citep{wu2011empirical,gayah2011effects, farhi2011traffic}, and loading and unloading processes \citep{gayah2011clockwise}. Efforts also have been devoted to deriving or approximating the shape of MFDs. Earlier in \citep{herman1979two}, a two-fluid model was introduced to town traffic by splitting vehicles into two groups: moving and stopped vehicles. Relations among flow, density, and speed were derived by assuming the average speed in an urban network is a function of the fraction of stopped vehicles at any given time. In \citep{daganzo2008analytical}, variational theory was used to approximate the MFD for any multi-block, signal-controlled streets without turning movements. This method was later extended in \citep{leclercq2013estimating} to provide an exact solution. In \citep{helbing2009derivation}, the MFD was derived by averaging the link-based fundamental diagrams for all links with the same parameter settings. In \citep{daganzo2011macroscopic-b}, a two-bin model was used to study the possible stationary states in a double-ring network, and multivaluedness was observed in the network flow-density relations. However, the aforementioned analytical studies have their limitations: route choice behaviors were not explicitly modeled in \citep{herman1979two,daganzo2008analytical,helbing2009derivation}, while signal settings were not explicitly modeled in \citep{herman1979two,daganzo2011macroscopic-b}. Furthermore, the definition of stationary states in signalized networks has not been clearly defined. 

In \citep{jin2013kinematic}, stationary states were analytically solved in an uninterrupted (i.e., unsignalized) double-ring network within the framework of kinematic wave models. It was found that infinitely many stable states with zero-speed shock waves (ZS) could arise at the same network density in the congested branch of MFDs.  Simulations using the cell transmission model (CTM) \citep{daganzo1994cell,daganzo1995cell} of a signalized double ring network unveiled periodic traffic patterns and these were defined as stationary states. However, due to the infinite-dimensional state variables introduced by the kinematic wave model, traffic statics and dynamics are very difficult to solve when traffic signals are present. To the best of our knowledge, the number of stationary states in a signalized double-ring network and their stability properties and relations to cycle lengths, retaining ratios, and initial densities are all still unclear.

\begin{figure}[!h]
\centering
\subfigure[A signalized intersection]{\includegraphics[scale=0.25]{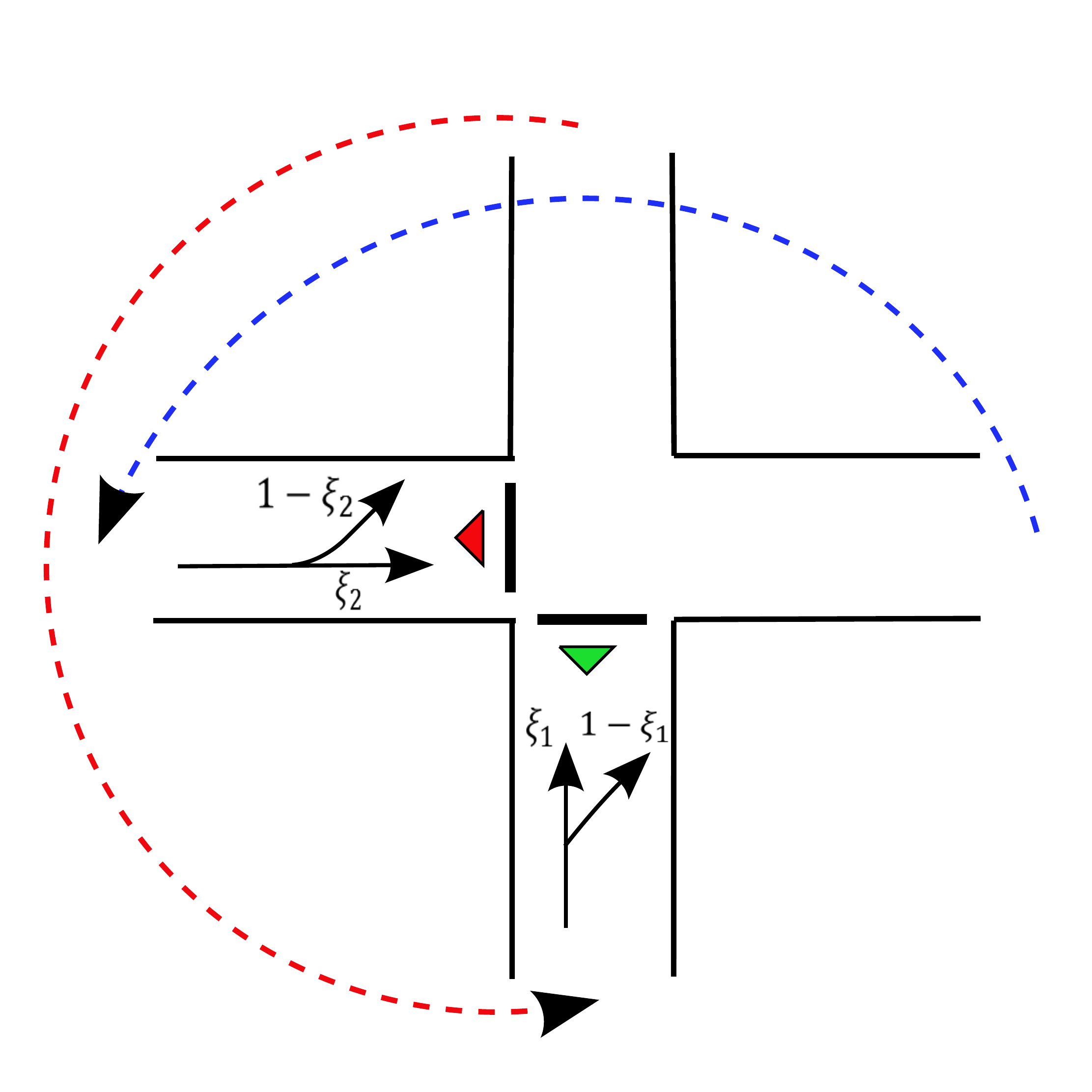} \label{fig:urban-network:intersection}}
\subfigure[An abstract network]{\includegraphics[scale=0.65]{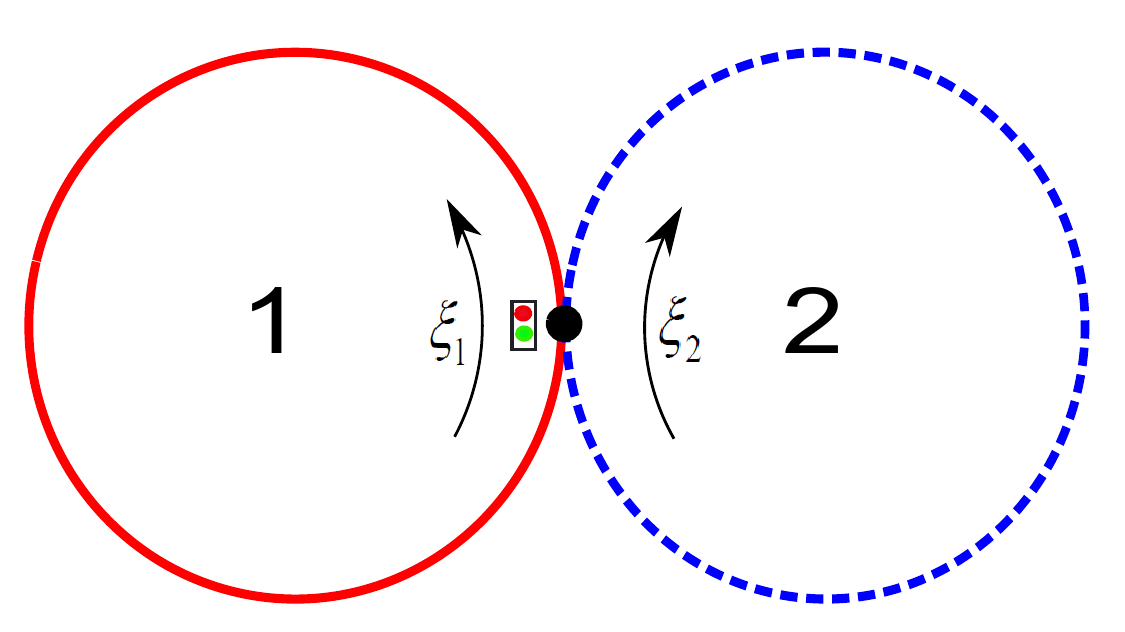} \label{fig:urban-network:double-ring}}
\caption{A signalized double-ring network. }
\label{fig:urban-network}
\end{figure}

This study aims to fill this gap by studying a simple signalized intersection shown in Figure \ref{fig:urban-network:intersection}. If exiting vehicles in the downstream are immediately added into the network from their corresponding upstream entrances, the signalized intersection is further changed into a signalized double-ring network with periodic boundary conditions, which is shown in Figure \ref{fig:urban-network:double-ring}. In the signalized double-ring network, traffic statics and dynamics are easier to solve because we only need to keep track of traffic states on the two rings.   
Furthermore, instead of using the kinematic wave model in \citep{jin2013kinematic}, we use a simplified approximation called the link queue model in \citep{jin2012link} to formulate the traffic dynamics in the signalized double-ring network. With the assumption of a triangular traffic flow fundamental diagram \citep{haberman1977mathematical}, the signalized double-ring network is studied as a switched affine system. The switching rule is governed by three sets of parameters: initial densities, route choice behaviors, and signal settings. Periodic density evolution orbits are found and defined as stationary states. A Poincar\'{e} map approach \citep{wiggins1990introduction,teschl2012ordinary} is used to study the properties of such stationary states. Poincar\'{e} maps were originally used to study the movements of celestial bodies and are a very important tool in analyzing periodic orbits: each periodic orbit corresponds to a fixed point on the Poincar\'{e} map, and its stability is directly related to the stability of the fixed point. In \citep{jin2013stability},  Poincar\'{e} maps have been applied to study the stability and bifurcation in traffic flow in diverge-merge and other networks within the framework of network kinematic wave theories. With short cycle lengths, closed-form Poincar\'{e} maps are derived for signalized networks, and stationary states and the corresponding stabilities are obtained by solving and analyzing the fixed points on the Poincar\'{e} maps. The impacts of retaining ratios and initial densities on the MFDs and the gridlock times are analyzed. With long cycle lengths, closed-form Poincar\'{e} maps are hard to obtain, and thus, the corresponding fixed points are numerically solved using the secant method \citep{epperson2014introduction}. 

The rest of this paper is organized as follows. In Section \ref{A switched affine system for a signalized double-ring network}, we formulate the traffic dynamics in a signalized double-ring network as a switched affine system when the link queue model \citep{jin2012link} and a triangular traffic flow fundamental diagram \citep{haberman1977mathematical} are used. In Section \ref{Periodic density evolution orbits and derivation of Poincare maps}, we derive the Poincar\'{e} maps from the density evolution orbits in the switched affine system. Stationary states and their stability properties are defined in terms of the fixed points and their stability properties on the Poincar\'{e} maps. In Section \ref{Stationary states and their stability properties with short cycle lengths}, we solve the stationary states and the corresponding stability properties with short cycle lengths. In Section \ref{Macroscopic fundamental diagrams and gridlock time analysis}, we discuss the impacts of retaining ratios and initial densities on the MFDs and the gridlock times.  In Section \ref{Numerical solutions to the Poincare maps with long cycle lengths}, we provide numerical solutions to the fixed points on the Poincar\'{e} maps when the cycle lengths are long. In Section \ref{Conclusions}, we draw our conclusions with some future research directions.

\section{A switched affine system for a signalized double-ring network} \label{A switched affine system for a signalized double-ring network}
\subsection{A link queue representation} \label{A link queue representation}
According to the link queue model (LQM) \citep{jin2012link}, vehicles in a link are treated as always being in a queue, and therefore, each link has only one state variable, the average link density. The left and the right rings in Figure \ref{fig:urban-network:double-ring} are denoted as rings 1 and 2 with average densities $k_1(t)$ and $k_2(t)$, respectively. We assume both rings have the same length $L$. Since it is a closed network, we have $k_1(t)+k_2(t)=2k$, where $k$ is the average network density. In addition, we assume both rings have the same location-and-time independent fundamental diagram, $q_a = Q(k_a)$, $a=1,2$, for $k_a\in[0,k_j]$, where $k_j$ is the jam density. Generally speaking, $Q(k_a)$ is a concave function with its capacity $C$ obtained at the critical density $k_{c}$, i.e., $C=Q(k_c)$. Then the local demand and supply are defined as
\begin{subequations}
\begin{align}
D_a(t)&=Q(\min\{k_a(t),k_c\})=\cas{{ll}Q(k_a(t)), & k_a(t)\in [0,k_c],\\ C, & k_a(t)\in (k_c,k_j],}  & a=1,2 \label{equ:demand:double-ring}\\
S_a(t)&=Q(\max\{k_a(t),k_c\})=\cas{{ll}C, & k_a(t)\in [0,k_c],\\ Q(k_a(t)), & k_a(t)\in (k_c,k_j],}  & a=1,2 \label{equ:supply:double-ring}
\end{align}
\label{equ:demand:supply:double-ring}
\end{subequations}

At the junction, vehicles in the two upstream approaches move alternately, and therefore, there are two phases in each cycle. Without loss of generality, we assign phase one to vehicles in ring 1 and phase two to vehicles in ring 2. The cycle length is $T$ with a lost time $\Delta$ for each phase. The green ratio is denoted as $\pi_1$ for ring 1 and $\pi_2$ for ring 2. We assume the yellow and all red period in each phase is the same as the lost time, and therefore, the effective green time is $\pi_1 T$ for ring 1 and $\pi_2 T$ for ring 2, and $(\pi_1+\pi_2)T=T-2\Delta$. Then the signal regulation can be described by the following two indicator functions:
\begin{subequations}
\begin{align}
\delta_1(t;T,\Delta,\pi_1)&=\cas{{ll}1, & t\in[nT,nT+\pi_1 T ),\\ 0, & \text{otherwise},}  & n\in\mathbb{N}_0 \label{equ:indicator function:ring1}\\
\delta_2(t;T,\Delta,\pi_1)&=\cas{{ll}1, & t\in[nT+\Delta+\pi_1 T,(n+1)T-\Delta ),\\0, & \text{otherwise},}  & n\in\mathbb{N}_0 \label{equ:indicator function:ring2}
\end{align}
\label{equ:indicator function}
\end{subequations}
where $\mathbb{N}_0=\{0,1,2,3,...\}$. According to Equation \refe{equ:indicator function}, we have three different combinations of $(\delta_1(t),\delta_2(t))$ within one cycle: (i) $(1,0)$ stands for the effective green time period in phase one; (ii) $(0,1)$ stands for the effective green time period in phase two; (iii) $(0,0)$ stands for the lost time period in either of the phases. 

Due to the signal regulations, the signalized $2\times 2$ junction is equivalent to two alternating diverging junctions, and the invariant first-in-first-out (FIFO) diverging model \citep{daganzo1995cell} is used. The retaining ratio is denoted as $\xi_1(t) \in(0,1)$ for ring 1 and $\xi_2(t)\in(0,1)$ for ring 2, which means the turning ratio is $1-\xi_1(t)$ from ring 1 to ring 2 and $1-\xi_2(t)$ from ring 2 to ring 1. Then the out-fluxes $g_1(t)$, $g_2(t)$, and the in-fluxes $f_1(t)$, $f_2(t)$ can be calculated as
\begin{subequations}
\begin{align}
g_1(t)&=\delta_1(t)\min\{D_1(t),\frac{S_1(t)}{\xi_1(t)},\frac{S_2(t)}{1-\xi_1(t)}\},\\
g_2(t)&=\delta_2(t)\min\{D_2(t),\frac{S_2(t)}{\xi_2(t)},\frac{S_1(t)}{1-\xi_2(t)}\},\\
f_1(t)&=g_1(t)\xi_1(t)+g_2(t)(1-\xi_2(t)),\\
f_2(t)&=g_1(t)(1-\xi_1(t))+g_2(t)\xi_2(t).
\end{align}
\label{equ:diverge:double-ring}
\end{subequations}

Since it is a closed network, we only need to focus on the traffic dynamics in one of the rings, e.g., ring 1. According to the traffic conservation in ring 1, the following equation holds:
\begin{equation}
\frac{d k_1(t)}{dt}=\frac{1}{L}(f_1(t)-g_1(t))=\frac{-(1-\xi_1(t))}{L}g_1(t)+\frac{(1-\xi_2(t))}{L}g_2(t).
\label{equ:link queue:double-ring}
\end{equation}
With Equations \refe{equ:demand:supply:double-ring} to \refe{equ:link queue:double-ring}, the system equation can be written as
\begin{equation}
\frac{d k_1(t)}{dt}=F(k_1(t);k,\delta_1(t),\delta_2(t),\xi_1(t),\xi_2(t)),
\label{equ:link queue:double-ring:system equation}
\end{equation}
which is closely related to the average network density $k$, the signal settings $(\delta_1(t),\delta_2(t))$, and the route choice behaviors $(\xi_1(t),\xi_2(t))$. Note that Equation \refe{equ:link queue:double-ring:system equation} is a nonlinear ordinary differential equation (ODE) with periodic forces. It is simpler than the kinematic wave model (Equation (3) in \citep{jin2013kinematic}), which is a partial differential equation (PDE), but still quite challenging to solve.  

\subsection{A switched affine system}\label{Formulation as a switched affine system}
For the remainder of this paper, we adopt the following triangular traffic flow fundamental diagram \citep{haberman1977mathematical}: 
\begin{equation}
Q(k_a)=\min\{v_fk_a, w_{a}(k_j-k_a)\}, \qquad a=1,2,
\label{equ:triangular FD}
\end{equation}
where $v_f$ and $w_{a}$ are the free-flow speed and the shock-wave speed, respectively. Then the signalized double-ring network can be further formulated as a switched affine system \citep{el2007poincare,sun2011stability}, and the traffic dynamics in Equation \refe{equ:link queue:double-ring:system equation} can be rewritten as
\begin{equation}
\frac{d k_1(t)}{dt}=A_{i}k_1(t)+B_{i}, \qquad i=1,2,...,N,
\label{equ:switched affine system}
\end{equation}
where $N$ is the number of possible combinations of $(A_i,B_i)$. The choice of $i$ (or $(A_i,B_i)$) is governed by three sets of parameters: the initial densities $(k_1(t),k)$, the signal settings $(\delta_1(t),\delta_2(t))$, and the route choice behaviors $(\xi_1(t),\xi_2(t))$. 

\begin{figure}[!h]
\centering
\subfigure[$(\delta_1(t),\delta_2(t))=(1,0)$]{\includegraphics[scale=0.23]{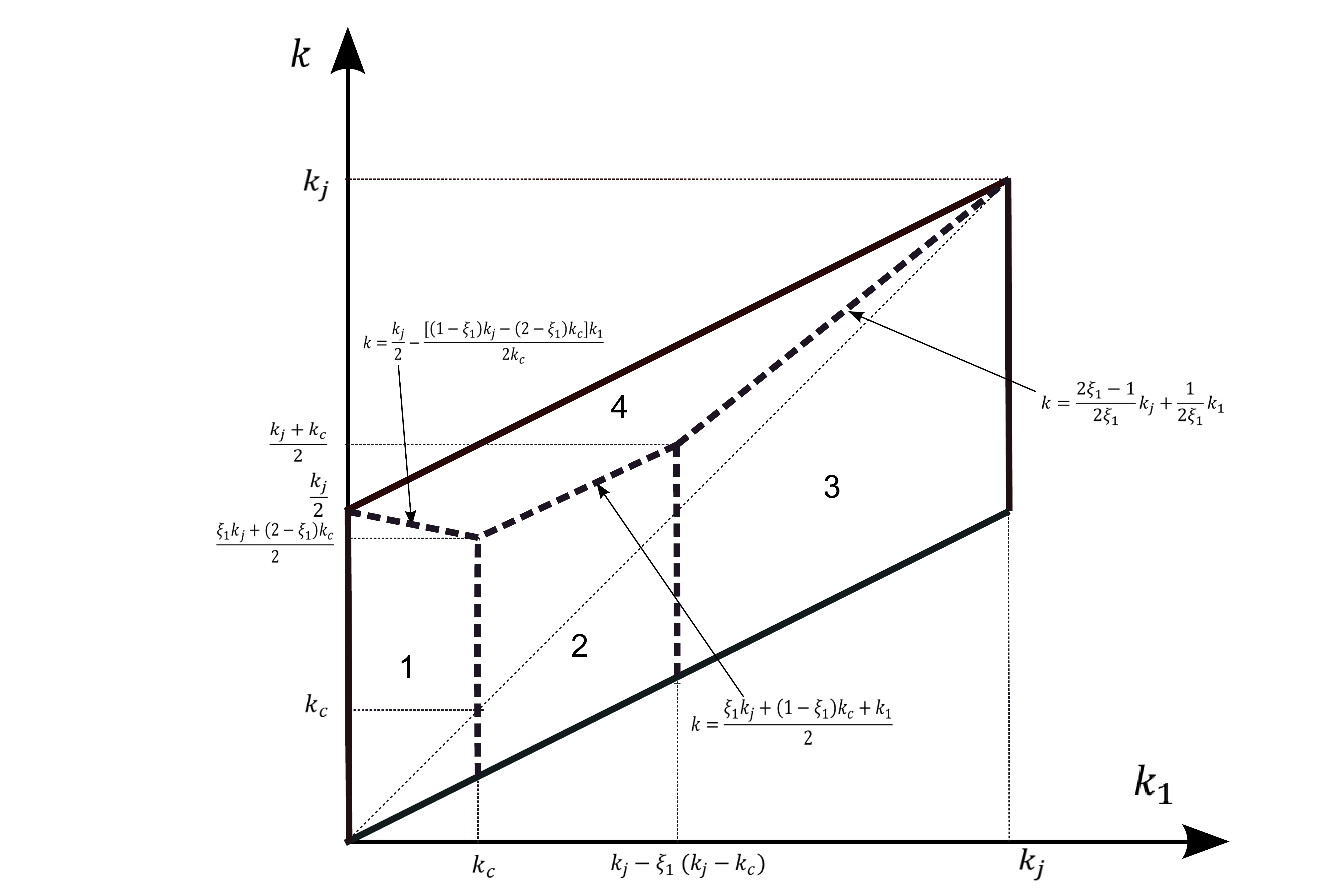}\label{fig:phase_diagram:delta1}}
\subfigure[$(\delta_1(t),\delta_2(t))=(0,1)$]{\includegraphics[scale=0.23]{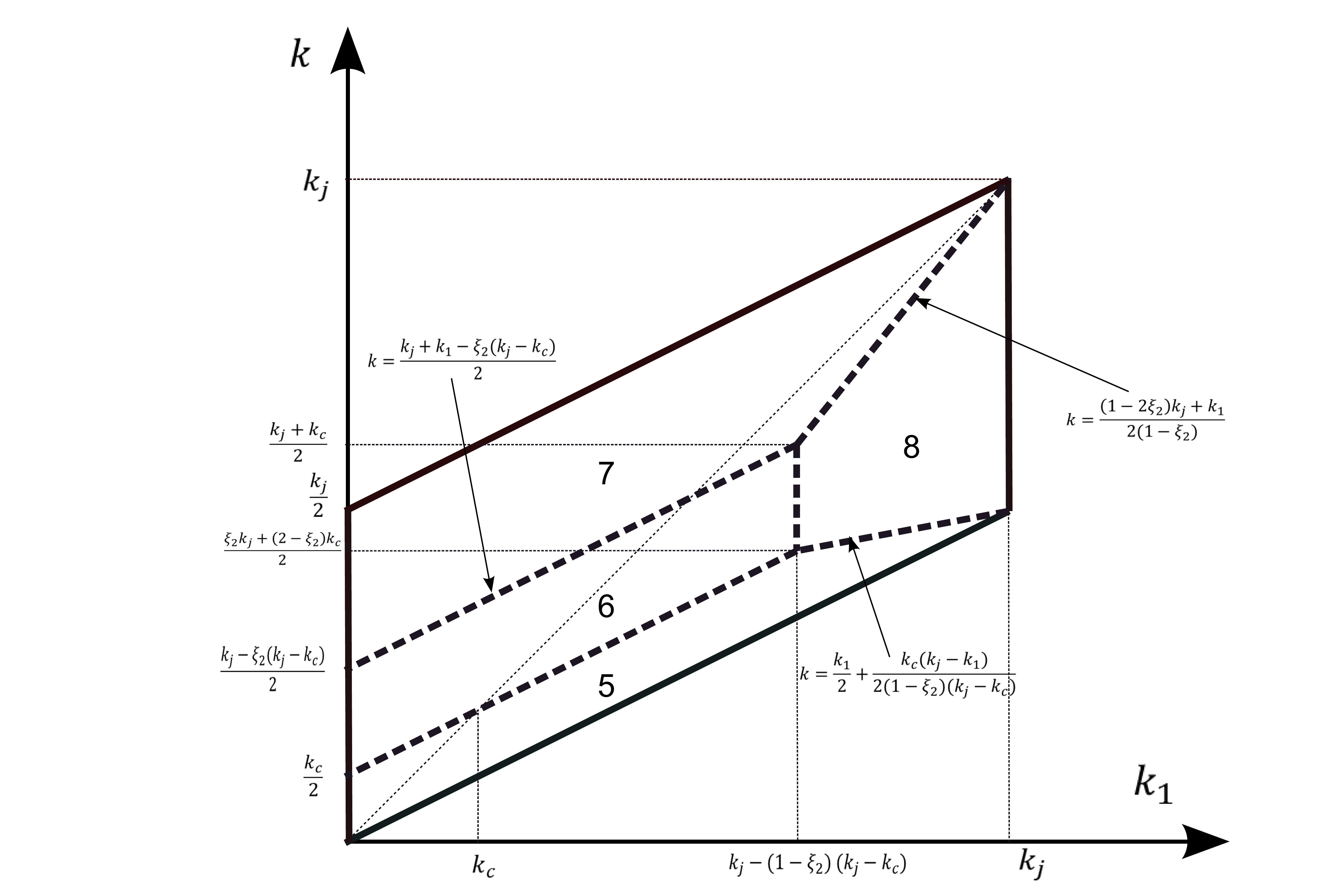}\label{fig:phase_diagram:delta2}}
\caption{Regions in the $(k_1,k)$ space.}
\label{fig:phase_diagram}
\end{figure}

Since densities in both rings vary during the effective green times, the $(k_1,k)$ space can be divided into different regions, in which the coefficients $(A_{i}, B_{i})$ remain the same. According to Equations \refe{equ:demand:supply:double-ring}, \refe{equ:diverge:double-ring}, and \refe{equ:triangular FD}, the $(k_1,k)$ space can be divided into four regions when $(\delta_1(t),\delta_2(t))=(1,0)$, which are provided in Figure \ref{fig:phase_diagram:delta1}. The meanings of the bold dashed lines in Figure \ref{fig:phase_diagram:delta1} are explained as follows:
\begin{itemize}
\item[(i)] the line $k=\frac{k_j}{2}-\frac{[(1-\xi_1)k_j-(2-\xi_1)k_c]k_1}{2k_c}$ stands for the case when $D_1=\frac{S_2}{1-\xi_1}<C<\frac{S_1}{\xi_1}$;
\item[(ii)] the line $k=\frac{\xi_1k_j+(1-\xi_1)k_c+k_1}{2}$ stands for the case when $D_1=\frac{S_2}{1-\xi_1}=C<\frac{S_1}{\xi_1}$;
\item[(iii)] the line $k=\frac{2\xi_1-1}{2\xi_1}k_j+\frac{1}{2\xi_1}k_1$ stands for the case when $\frac{S_1}{\xi_1}=\frac{S_2}{1-\xi_1}<D_1=C$;
\item[(iv)]  the line $k_1=k_j-\xi_1(k_j-k_c)$ stands for the case when $D_1=\frac{S_1}{\xi_1}=C<\frac{S_2}{1-\xi_1}$.
\end{itemize}
Similarly, the $(k_1,k)$ space can be divided into another four regions when $(\delta_1(t),\delta_2(t))=(0,1)$, which are provided in Figure \ref{fig:phase_diagram:delta2}. The meanings of the bold dashed lines in Figure \ref{fig:phase_diagram:delta2} are provided below:
\begin{itemize}
\item[(i)] the line $k=\frac{k_j+k_1-\xi_2(k_j-k_c)}{2}$ stands for the case when $D_2=\frac{S_2}{\xi_2}=C<\frac{S_1}{1-\xi_2}$;
\item[(ii)] the line $k=\frac{(1-2\xi_2)k_j+k_1}{2(1-\xi_2)}$ stands for the case when $\frac{S_2}{\xi_2}=\frac{S_1}{1-\xi_2}<D_2=C$;
\item[(iii)] the line $k=\frac{k_1}{2}+\frac{k_c(k_j-k_1)}{2(1-\xi_2)(k_j-k_c)}$ stands for the case when $D_2=\frac{S_1}{1-\xi_2}<C<\frac{S_2}{\xi_2}$;
\item[(iv)]  the line $k_1=k_j-(1-\xi_2)(k_j-k_c)$ stands for the case when $D_2=\frac{S_1}{(1-\xi_2)}=C<\frac{S_2}{\xi_2}$.
\end{itemize}
During the lost time periods (i.e., $(\delta_1(t),\delta_2(t))=(0,0)$) in phases one and two, densities in both rings remain unchanged, and therefore, we consider them as two different regions with the same coefficients in Equation \refe{equ:switched affine system} (i.e., $(A_{i}, B_{i})=(0,0)$). Thus, we have $N=10$. The possible values of $A_{i}$ and $B_{i}$ and the corresponding conditions are provided in Table \ref{table:possible coefficients}. 
\begin{table}[!h]\footnotesize
\centering
\caption{ Possible values of $A_{i}$ and $B_{i}$ and the corresponding conditions }
\begin{tabular}[c]{|c|c|c|c|c|}
\hline
{Region ($i$)}&{$A_{i}$}&{$B_{i}$}&{Conditions}\\
\hline
{1}&{$-\gamma_1$}&{0}&{$(\delta_1(t),\delta_2(t))=(1,0)$, $0<k_1<k_c$, $\frac{k_1}{2}\leq k\leq \frac{k_j}{2}-\frac{((1-\xi_1)k_j-(2-\xi_1)k_c)k_1}{2k_c}$}\\
\hline
{2}&{0}&{$-\gamma_1k_c$}&{$(\delta_1(t),\delta_2(t))=(1,0)$, $k_c\leq k_1<k_j-\xi_1(k_j-k_c)$, $\frac{k_1}{2}\leq k\leq\frac{\xi_1k_j+(1-\xi_1)k_c+k_1}{2}$}\\
\hline
{3}&{$\gamma_2$}&{$-\gamma_2 k_j$}&{$(\delta_1(t),\delta_2(t))=(1,0)$, $k_j-\xi_1(k_j-k_c)\leq k_1\leq k_j$, $\frac{k_1}{2}\leq k\leq\frac{2\xi_1 -1}{2\xi_1}k_j+\frac{1}{2\xi_1}k_1$}\\
\hline
\multirow{2}{*}{4}&\multirow{2}{*}{$-\gamma_3$}&\multirow{2}{*}{$-\gamma_3(k_j-2k)$}&{$(\delta_1(t),\delta_2(t))=(1,0)$, $0< k_1< k_j$}\\
{}&{}&{}&{$\max\{\frac{k_j}{2}-\frac{[(1-\xi_1)k_j-(2-\xi_1)k_c]k_1}{2k_c},\frac{\xi_1k_j+(1-\xi_1)k_c+k_1}{2},\frac{2\xi_1 -1}{2\xi_1}k_j+\frac{1}{2\xi_1}k_1 \} < k\leq\frac{k_j+k_1}{2}$}\\
\hline
{5}&{-$\gamma_4$}&{$2k\gamma_4$}&{$(\delta_1(t),\delta_2(t))=(0,1)$, $0\leq k_1\leq k_j$, $\frac{k_1}{2}\leq k\leq\min\{\frac{k_1+k_c}{2},\frac{k_1}{2}+\frac{k_c(k_j-k_1)}{2(1-\xi_2)(k_j-k_c)}\}$}\\
\hline
{6}&{0}&{$\gamma_4k_c$}&{$(\delta_1(t),\delta_2(t))=(0,1)$, $0\leq k_1\leq k_j-(1-\xi_2)(k_j-k_c)$, $\frac{k_1+k_c}{2}< k\leq \frac{k_j+k_1-\xi_2(k_j-k_c)}{2}$}\\
\hline
{7}&{$\gamma_5$}&{$\gamma_5(k_j-2k)$}&{$(\delta_1(t),\delta_2(t))=(0,1)$, $0\leq k_1\leq k_j$, $\max\{\frac{k_j+k_1-\xi_2(k_j-k_c)}{2},\frac{(1-2\xi_2)k_j+k_1}{2(1-\xi_2)}\}< k \leq \frac{k_1+k_j}{2}$}\\
\hline
\multirow{2}{*}{8}&\multirow{2}{*}{$-\gamma_3$}&\multirow{2}{*}{$\gamma_3k_j$}&{$(\delta_1(t),\delta_2(t))=(0,1)$, $k_j-(1-\xi_2)(k_j-k_c) < k_1\leq k_j$}\\
{}&{}&{}&{$\frac{k_1}{2}+\frac{k_c(k_j-k_1)}{2(1-\xi_2)(k_j-k_c)}< k < \frac{(1-2\xi_2)k_j+k_1}{2(1-\xi_2)}$}\\
\hline
{$9$}&{$0$}&{$0$}&{$(\delta_1(t),\delta_2(t))=(0,0)$, and transition from $(\delta_1(t),\delta_2(t))=(1,0)$ }\\
\hline
{$10$}&{$0$}&{$0$}&{$(\delta_1(t),\delta_2(t))=(0,0)$, and transition from $(\delta_1(t),\delta_2(t))=(0,1)$ }\\
\hline
\end{tabular}
\label{table:possible coefficients}
\begin{flushleft}{where $\gamma_1=\frac{(1-\xi_1)v_f}{L}$,  $\gamma_2=\frac{(1-\xi_1)v_fk_c}{L\xi_1(k_j-k_c)}$, $\gamma_3=\frac{v_fk_c}{L(k_j-k_c)}$, $\gamma_4=\frac{(1-\xi_2)v_f}{L}$, and $\gamma_5=\frac{(1-\xi_2)v_fk_c}{L\xi_2(k_j-k_c)}$.}
\end{flushleft}
\end{table}

\section{Periodic density evolution orbits and derivation of Poincar\'{e} maps} \label{Periodic density evolution orbits and derivation of Poincare maps}
\subsection{Periodic density evolution orbits} \label{Periodic density evolution orbits}
\begin{figure}[!h]
\centering
\subfigure[In the $(k_1,k)$ space]{\includegraphics[scale=0.18]{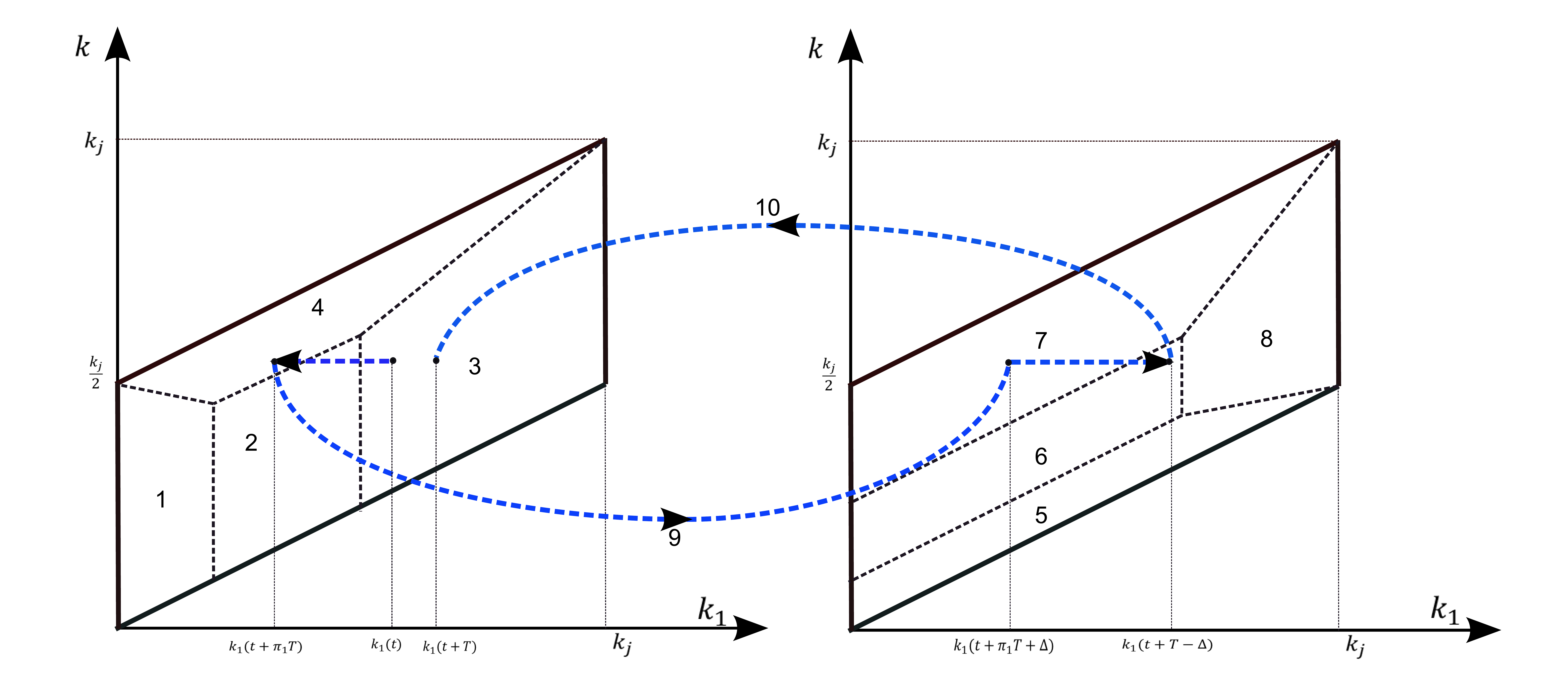}\label{fig:phase_periodic_maps}}
\subfigure[In the mapping of $k_1$]{\includegraphics[scale=0.2]{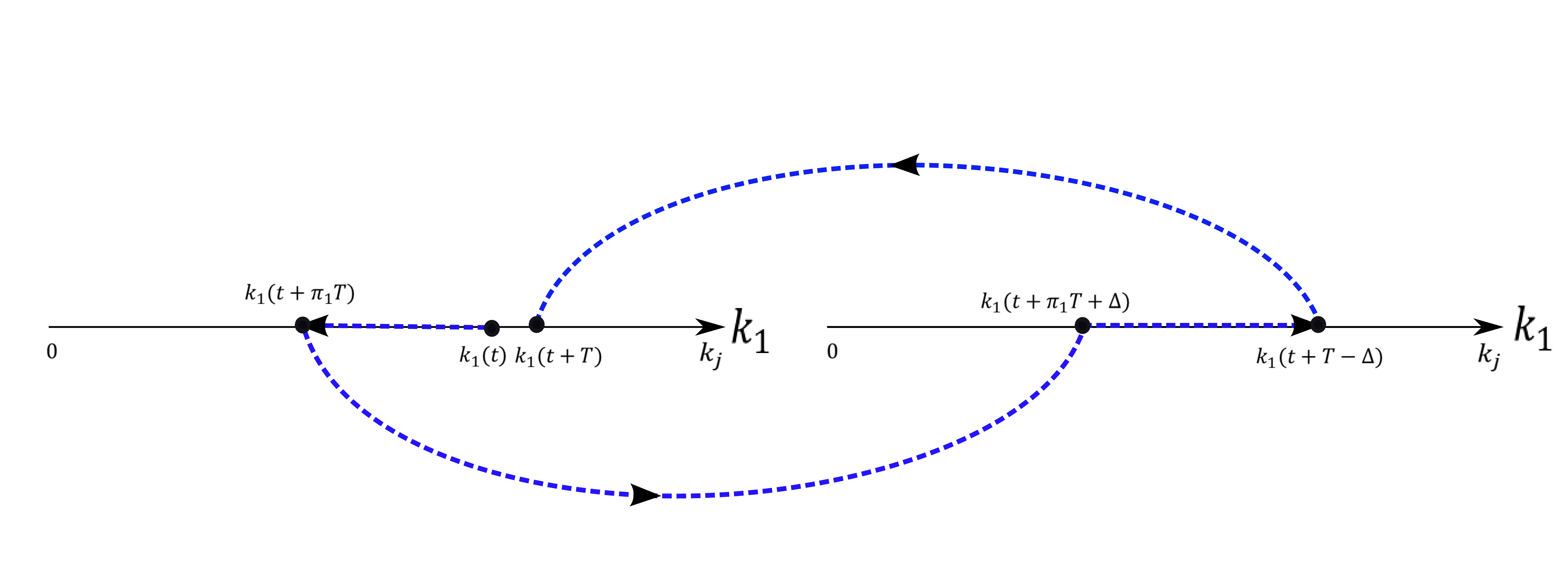}\label{fig:phase_poincare_maps}}
\caption{The density evolution orbit within one cycle.}
\label{fig:phase_periodic_poincare_maps}
\end{figure}

Due to periodic signal regulations, the switched affine system periodically visits the $(k_1,k)$ space in Figure \ref{fig:phase_diagram}, and therefore, density evolution orbits are formed. The circulating period is fixed, which is the cycle length T. In Figure \ref{fig:phase_periodic_maps}, we provide the density evolution orbit within one cycle in the $(k_1,k)$ space. And also, since it is a closed network (i.e., $k$ is fixed once it is given), we can map the density evolution orbit to the $k_1$ axis, which is shown in Figure \ref{fig:phase_poincare_maps}. The evolution process can be described as follows:
\begin{itemize}
\item [(1)] At the beginning of the cycle at time $t$, the initial densities are denoted as $(k_1(t),k)$. During $[t,t+\pi_1T)$, we have $(\delta_1,\delta_2)=(1,0)$. The density $k_1$ decreases as time elapses since vehicles in ring 1 can diverge to either of the rings while vehicles in ring 2 have to wait. The density evolution may cross multiple regions during this time period, and $i$ can take multiple values from 1 to 4 in Table \ref{table:possible coefficients}.  
\item [(2)] At time $t+\pi_1T$, the system switches to region 9. During $[t+\pi_1T, t+\pi_1 T+\Delta)$, we have $(\delta_1,\delta_2)=(0,0)$. The densities in both rings remain the same since no vehicles can move to the downstream links.  
\item [(3)] At time $t+\pi_1 T+\Delta$, the system switches to phase two. During $[t+\pi_1 T+\Delta, t+T-\Delta)$, we have $(\delta_1,\delta_2)=(0,1)$. The density $k_1$ increases as time elapses since vehicles in ring 2 are allowed to diverge to either of the rings while vehicles in ring 1 have to wait. Similarly, the density evolution may cross multiple regions during this time period, and $i$ can take multiple values from 5 to 8 in Table \ref{table:possible coefficients}. 
\item [(4)] At time $t+T-\Delta$, the system switches to region 10. During $[t+T-\Delta, t+T)$,  $(\delta_1,\delta_2)=(0,0)$ and the density in ring 1 remains the same.
\item [(5)] At $t+T$, the system switches back to Step (1) with $(\delta_1,\delta_2)=(1,0)$. The densities become $(k_1(t+T),k)$ and serve as new initial densities in the next cycle.
\end{itemize}

Based on the above description, it is possible for the signalized double-ring network to have periodic density evolution orbits if $k_1(t+T)=k_1(t)$. Such periodic patterns were observed in the simulations  and defined as stationary states in \citep{jin2013kinematic}. To calculate the asymptotic average network flow-rate, the following equation was used:
\begin{equation}
q(t)=\frac{\hat{g}_1(t)+\hat{g}_2(t)}{2}=\frac{\int_{s=t-T}^{t}g_1(s)ds+\int_{s=t-T}^{t}g_2(s)ds}{2T}.
\label{equ:average flow-rate}
\end{equation} 
When the system reaches a stationary state, $q(t)$ becomes a constant value, which depends on the average network density; i.e., $q(t)|_{t\rightarrow +\infty}=q(k)$. Then the relation between $q(k)$ and $k$ is the MFD for the signalized double-ring network. 

\subsection{Derivation of Poincar\'{e} maps} \label{Derivation of Poincare maps}
In the following sections we study the static and dynamic properties associated with the density evolution orbits in terms of Poincar\'{e} maps. According to the discussion in Section \ref{Periodic density evolution orbits}, the density evolution orbit for each cycle is combined with the following four local maps, $P_i$, $i=1,2,3,4$:
\begin{subequations}
\begin{align}
k_1(t+\pi_1 T) &=P_1k_1(t), \\
k_1(t+\pi_1 T+\Delta) &=P_2k_1(t+\pi_1 T), \\
k_1(t+T-\Delta ) &=P_3k_1(t+\pi_1 T+\Delta), \\
k_1(t+T) &=P_4k_1(t+T-\Delta ).
\end{align}
\label{equ:poincare_map:local maps}
\end{subequations}
If we define the Poincar\'{e} section as the time when the system first visits the $(k_1,k)$ space in Figure \ref{fig:phase_periodic_maps} in each cycle, the Poincar\'{e} map can be derived as the composition of the four local mappings.
\begin{equation}
k_1(t+T)=Pk_1(t)=P_4\circ P_3\circ P_2\circ P_1k_1(t).
\label{equ:poincare_map:equation}
\end{equation}
The Poincar\'{e} map in Equation \refe{equ:poincare_map:equation} is well defined and can be analytically derived if initial densities, retaining ratios, and signal settings are given. Then we can have the following definition of stationary states:
\begin{definition}[Stationary states]\label{definition:stationary states: fixed points}
The signalized double-ring network is in a stationary state when there exists a fixed point $k_1^{*}$ on the Poincar\'{e} map $P k_1(t)$ satisfying $k_1^{*}=P k_1^{*}$.
\end{definition}

Due to the existence of noise in transportation networks, it is necessary to understand the stability properties of stationary states. Suppose at the beginning of one cycle, there exists a perturbation that drives the density in ring 1 from $k_1^{*}$ to  $k_1(t)=k_1^{*}+\epsilon(t)$, where $\epsilon(t)$ is a small error term. Then after n cycles, the density in ring 1 and the error term become
\begin{subequations}
\begin{align}
k_1(t+nT)&=P^nk_1(t), &\qquad n\in\mathbb{N}_0\\
\epsilon(t+nT)&=P^nk_1(t)-k_1^{*}=(\partial P)^n \epsilon (t), &\qquad n\in\mathbb{N}_0\label{equ:poincare_map:n cycle:error term}
\end{align}
\label{equ:poincare_map:n cycle}
\end{subequations} 
Following \citep{la1976stability,teschl2012ordinary}, we have the following definition of stability for the fixed points on the the Poincar\'{e} maps:
\begin{definition}[Local stability]\label{definition:local stability}
The fixed point $k_1^{*}$ is: (i) Lyapunov stable if for any given $\beta>0$, there exists an $\omega\leq \beta$ such that $|\epsilon(t+nT)|< \beta$ for any $|\epsilon(t)|<\omega$ and $n\in\mathbb{N}_0$; (ii) unstable otherwise. Furthermore, $k_1^{*}$ is asymptotically stable if it is stable and $\epsilon(t+nT) \rightarrow 0$ as $n\rightarrow + \infty$.
\end{definition}

\section{Stationary states and their stability properties with short cycle lengths} \label{Stationary states and their stability properties with short cycle lengths}
For simplicity, we consider the following homogeneous settings in the rest of the paper: both rings have the same effective green times,  i.e., $\pi_2=\pi_1=\pi$, and the same time-independent retaining ratios, i.e., $\xi_1=\xi_2=\xi$. 

\subsection{Possible stationary states}\label{Possible stationary states}
To derive closed-form Poincar\'{e} maps, we require the cycle length to be short enough so that it guarantees the density evolution orbit doesn't cross multiple regions during each effective green time period, i.e., the density evolution orbit stays only within region $i\in \{1,2,3,4\}$ when $(\delta_1,\delta_2)=(1,0)$, and only in region $j\in\{5,6,7,8\}$ when $(\delta_1,\delta_2)=(0,1)$. Under different retaining ratios, e.g., $\xi>0.5$ and $\xi < 0.5$, we have 11 different combinations of regions in one cycle, which are provided in Table \ref{table:density evolution signs for 11 regions}. Here, we define the sum of the coefficients $(A_{i},B_{i})$ and $(A_{j},B_{j})$ as
\begin{equation}
\lambda(k_1,k)=A_{j}k_1+B_{j}+A_{i}k_1+B_{i}. \label{equ:indicator function:phase diagram}
\end{equation}
When the network is not gridlocked, $k_1$ will keep increasing or decreasing if $\lambda(k_1,k)$ is always greater or smaller than zero, which means it is impossible to have stationary states with the combination of regions $(i,j)$. But if $\lambda(k_1,k)$ can take both positive and negative values or is always zero, it is possible to have stationary states inside the combination of regions $(i,j)$. In Table \ref{table:density evolution signs for 11 regions}, we provide possible values of $\lambda(k_1,k)$ and the corresponding conditions for different combinations of regions under different retaining ratios.
\begin{table}[!h]\small
\centering
\caption{ Possible values of $\lambda(k_1,k)$ in the 11 combinations of regions }
\begin{tabular}[c]{|c|c|c|c|c|c|}
\hline
\multicolumn{3}{|c|}{$0.5<\xi<1$}&\multicolumn{3}{c|}{$0<\xi < 0.5$}\\
\hline
{Regions}&{$\lambda(k_1,k)$}&{Condition}&{Regions}&{$\lambda(k_1,k)$}&{Condition}\\
\hline
{}&{$>0$}&{$k_1<k$}&{}&{$>0$}&{$k_1<k$}\\
\cline{2-3}\cline{5-6}
{$(1,5)$}&{$=0$}&{$k_1=k$}&{$(1,5)$}&{$=0$}&{$k_1=k$}\\
\cline{2-3}\cline{5-6}
{}&{$<0$}&{$k_1>k$}&{}&{$<0$}&{$k_1>k$}\\
\hline
{$(1,6)$}&{$>0$}&{}&{$(1,6)$}&{$>0$}&{}\\
\hline
{}&{$>0$}&{$\gamma_1<\gamma_5$ or}&{}&\multirow{3}{*}{}&{}\\
{}&{}&{$\gamma_1>\gamma_5$ and $k_1<\frac{\gamma_5(k_j-2k)}{\gamma_1-\gamma_5}$}&{}&{}&{}\\
\cline{2-3}
{$(1,7)$}&{$=0$}&{$\gamma_1>\gamma_5$ and $k_1=\frac{\gamma_5(k_j-2k)}{\gamma_1-\gamma_5}$}&{$(1,7)$}&{$>0$}&{}\\
\cline{2-3}
{}&{$<0$}&{$\gamma_1>\gamma_5$ and $k_1>\frac{\gamma_5(k_j-2k)}{\gamma_1-\gamma_5}$}&{}&{}&{}\\
\hline
{$(2,5)$}&{$<0$}&{}&{$(2,5)$}&{$<0$}&{}\\
\hline
{$(2,6)$}&{$=0$}&{}&{$(2,6)$}&{$=0$}&{}\\
\hline
{$(2,7)$}&{$<0$}&{}&{$(4,6)$}&{$>0$}&{}\\
\hline
{$(4,7)$}&{$<0$}&{}&{$(4,7)$}&{$>0$}&{}\\
\hline
{}&{$>0$}&{$\gamma_2<\gamma_4$ and $k_1<\frac{2k\frac{\xi(k_j-k_c)}{k_c}-k_j}{\frac{\xi(k_j-k_c)}{k_c}-1}$}&{}&\multirow{3}{*}{}&{}\\
\cline{2-3}
{$(3,5)$}&{$=0$}&{$\gamma_2<\gamma_4$ and $k_1=\frac{2k\frac{\xi(k_j-k_c)}{k_c}-k_j}{\frac{\xi(k_j-k_c)}{k_c}-1}$}&{$(3,5)$}&{$<0$}&{}\\
\cline{2-3}
{}&{$<0$}&{$\gamma_2<\gamma_4$ and $k_1>\frac{2k\frac{\xi(k_j-k_c)}{k_c}-k_j}{\frac{\xi(k_j-k_c)}{k_c}-1}$}&{}&{}&{}\\
{}&{}&{or $\gamma_2>\gamma_4$}&{}&{}&{}\\
\hline
{$(3,6)$}&{$>0$}&{}&{$(2,8)$}&{$<0$}&{}\\
\hline
{}&{$>0$}&{$k_1>k$}&{}&{$>0$}&{$k_1<k$}\\
\cline{2-3}\cline{5-6}
{$(3,7)$}&{$=0$}&{$k_1=k$}&{$(4,8)$}&{$=0$}&{$k_1=k$}\\
\cline{2-3}\cline{5-6}
{}&{$<0$}&{$k_1<k$}&{}&{$<0$}&{$k_1>k$}\\
\hline
{$(3,8)$}&{$>0$}&{}&{$(3,8)$}&{$<0$}&{}\\
\hline
\end{tabular}
\label{table:density evolution signs for 11 regions}
\end{table}
Specifically, when $\xi=0.5$, we only have the following combinations of regions $(i,j)$: $(1,5)$, $(1,6)$, $(1,7)$, $(2,5)$, $(2,6)$, $(4,7)$, $(3,5)$, and $(3,8)$. The values of $\lambda(k_1,k)$ and the corresponding conditions in regions $(1,5)$, $(1,6)$, $(1,7)$, $(2,5)$, $(2,6)$, and $(3,5)$ are the same as those with $\xi<0.5$ in Table \ref{table:density evolution signs for 11 regions}. However, the values of $\lambda(k_1,k)$ in regions $(4,7)$ and $(3,8)$ are always zero with $\xi=0.5$. Therefore, we have the following lemma:

\begin{lemma}[Possible regions having stationary states]
When the network is not gridlocked, it is only possible for the following combinations of regions $(i,j)$ to have stationary states:
\begin{itemize}
\item [(i)] $(1,5)$, $(1,7)$, $(2,6)$, $(3,5)$, and $(3,7)$ for $0.5<\xi<1$;
\item [(ii)] $(1,5)$, $(2,6)$, and $(4,8)$ for $0<\xi < 0.5$;
\item [(iii)] $(1,5)$, $(2,6)$, $(4,7)$, and $(3,8)$ for $\xi= 0.5$.
\end{itemize} 
However, when the network is gridlocked with one or both rings jammed, it is only possible for the combinations of regions  $(4,7)$ and $(3,8)$ to have stationary states with any $\xi\in(0,1)$.
\end{lemma}
The proof is simple and thus omitted here. Note that when the network is gridlocked with one or both rings jammed, the gridlock states are stationary states since they satisfy Definition \ref{definition:stationary states: fixed points}. It is easy to identify that these stationary states are inside the combinations of regions $(4,7)$ and $(3,8)$.

After identifying the possible combinations of regions having stationary states, we can analytically derive the closed-form Poincar\'{e} maps. The corresponding fixed points are the stationary states we are interested in. Therefore, we have the following theorem:
\begin{theorem}[Possible stationary states]\label{theorem:possible stationary states}
According to Equation \refe{equ:switched affine system} and Table \ref{table:possible coefficients}, the Poincar\'{e} maps and the corresponding fixed points in the possible regions having stationary states are derived and provided in Table \ref{table:possible poincare maps}.
\begin{table}[!h]\small
\centering
\caption{ Poincar\'{e} maps and fixed points in the possible regions having stationary states }
\begin{tabular}[c]{|c|c|c|}
\hline
{Regions}&{Poincar\'{e} map $Pk_1(t)$ with $0.5<\xi<1$}&{Fixed points $k_1^{*}$}\\
\hline
{$(1,5)$}&{$2k(1-e^{-\gamma_1\pi T})+k_1(t)e^{-2\gamma_1\pi T}$}&{$\frac{2k}{1+e^{-\gamma_1\pi T}}$}\\
\hline
{$(1,7)$}&{$(k_j-2k)(e^{\gamma_5\pi T}-1)+k_1(t)e^{(\gamma_5-\gamma_1)\pi T}$}&{$\frac{(k_j-2k)(e^{\gamma_5\pi T}-1)}{1-e^{(\gamma_5-\gamma_1)\pi T}}$}\\
\hline
{$(2,6)$}&{$k_1(t)$}&{$k_1(t)$}\\
\hline
{$(4,7)$}&{$(k_j-2k)(e^{(\gamma_5-\gamma_3)\pi T}-1)+k_1(t)e^{(\gamma_5-\gamma_3)\pi T}$}&{$2k-k_j$}\\
\hline
{$(3,5)$}&{$k_j(1-e^{\gamma_2\pi T})e^{-\gamma_4\pi T}+2k(1-e^{-\gamma_4\pi T})+k_1(t)e^{(\gamma_2-\gamma_4)\pi T}$}&{$\frac{2k(1-e^{-\gamma_4\pi T})-k_j(e^{\gamma_2\pi T}-1)e^{-\gamma_4\pi T}}{1-e^{(\gamma_2-\gamma_4)\pi T}}$}\\
\hline
{$(3,7)$}&{$k_j(2e^{\gamma_2\pi T}-e^{2\gamma_2\pi T}-1)-2k(e^{\gamma_2\pi T}-1)+k_1(t)e^{2\gamma_2\pi T}$}&{$\frac{2k+k_j(e^{\gamma_2\pi T}-1)}{e^{\gamma_2\pi T}+1}$}\\
\hline
{$(3,8)$}&{$k_j(1-e^{(\gamma_2-\gamma_3)\pi T})+k_1(t)e^{(\gamma_2-\gamma_3)\pi T}$}&{$k_j$}\\
\hline
\hline
{Regions}&{Poincar\'{e} map $Pk_1(t)$ with $0<\xi <0.5$}&{Fixed points $k_1^{*}$}\\
\hline
{$(1,5)$}&{$2k(1-e^{-\gamma_1\pi T})+k_1(t)e^{-2\gamma_1\pi T}$}&{$\frac{2k}{1+e^{-\gamma_1\pi T}}$}\\
\hline
{$(2,6)$}&{$k_1(t)$}&{$k_1(t)$}\\
\hline
{$(4,7)$}&{$(k_j-2k)(e^{(\gamma_5-\gamma_3)\pi T}-1)+k_1(t)e^{(\gamma_5-\gamma_3)\pi T}$}&{$2k-k_j$}\\
\hline
{$(4,8)$}&{$k_j(e^{-2\gamma_3\pi T}-2e^{-\gamma_3\pi T}+1)-2k(e^{-2\gamma_3\pi T}-e^{-\gamma_3\pi T})$}&{$\frac{k_j(1-e^{-\gamma_3\pi T})+2k e^{-\gamma_3\pi T}}{1+e^{-\gamma_3\pi T}}$}\\
{}&{$+k_1(t)e^{-2\gamma_3\pi T}$}&{}\\
\hline
{$(3,8)$}&{$k_j(1-e^{(\gamma_2-\gamma_3)\pi T})+k_1(t)e^{(\gamma_2-\gamma_3)\pi T}$}&{$k_j$}\\
\hline
\hline
{Regions}&{Poincar\'{e} map $Pk_1(t)$ with $\xi= 0.5$}&{Fixed points $k_1^{*}$}\\
\hline
{$(1,5)$}&{$2k(1-e^{-\gamma_1\pi T})+k_1(t)e^{-2\gamma_1\pi T}$}&{$\frac{2k}{1+e^{-\gamma_1\pi T}}$}\\
\hline
{$(2,6)$}&{$k_1(t)$}&{$k_1(t)$}\\
\hline
{$(4,7)$}&{$k_1(t)$}&{$k_1(t)$}\\
\hline
{$(3,8)$}&{$k_1(t)$}&{$k_1(t)$}\\
\hline
\end{tabular}
\label{table:possible poincare maps}
\end{table}
\end{theorem}

\subsection{Stability properties of the stationary states} \label{Stability properties of the stationary states}
In this subsection, we analyze the stability properties of the stationary states provided in Table \ref{table:possible poincare maps}, particularly the gridlock states. According to Equation \refe{equ:poincare_map:n cycle:error term}, the error term changes to $\epsilon(t+T)=\partial P \epsilon (t)$ after one cycle. In Table \ref{table:possible poincare maps:perturbation}, we provide the error term $\epsilon(t+T)$ in different combinations of regions under different retaining ratios. According to Definition \ref{definition:local stability}, we have the following theorem.

\begin{table}[!h]
\centering
\caption{ Changes in the error term after one cycle}
\begin{tabular}[c]{|c|c|}
\hline
{Regions}&{$\epsilon(t+T)$ with $0.5<\xi<1$}\\
\hline
{$(1,5)$}&{$\epsilon(t)e^{-2\gamma_1\pi T}$}\\
\hline
{$(1,7)$}&{$\epsilon(t)e^{(\gamma_5-\gamma_1)\pi T}$}\\
\hline
{$(2,6)$}&{$\epsilon(t)$}\\
\hline
{$(4,7)$}&{$\epsilon(t)e^{(\gamma_5-\gamma_3)\pi T}$}\\
\hline
{$(3,5)$}&{$\epsilon(t)e^{(\gamma_2-\gamma_4)\pi T}$}\\
\hline
{$(3,7)$}&{$\epsilon(t)e^{2\gamma_2\pi T}$}\\
\hline
{$(3,8)$}&{$\epsilon(t)e^{(\gamma_2-\gamma_3)\pi T}$}\\
\hline
\hline
{Regions}&{$\epsilon(t+T)$ with $0<\xi< 0.5$}\\
\hline
{$(1,5)$}&{$\epsilon(t)e^{-2\gamma_1\pi T}$}\\
\hline
{$(2,6)$}&{$\epsilon(t)$}\\
\hline
{$(4,7)$}&{$\epsilon(t)e^{(\gamma_5-\gamma_3)\pi T}$}\\
\hline
{$(4,8)$}&{$\epsilon(t)e^{-2\gamma_3\pi T}$}\\
\hline
{$(3,8)$}&{$\epsilon(t)e^{(\gamma_2-\gamma_3)\pi T}$}\\
\hline
\hline
{Regions}&{$\epsilon(t+T)$ with $\xi= 0.5$}\\
\hline
{$(1,5)$}&{$\epsilon(t)e^{-2\gamma_1\pi T}$}\\
\hline
{$(2,6)$}&{$\epsilon(t)$}\\
\hline
{$(4,7)$}&{$\epsilon(t)$}\\
\hline
{$(3,8)$}&{$\epsilon(t)$}\\
\hline
\end{tabular}
\label{table:possible poincare maps:perturbation}
\end{table}

\begin{theorem}[Stability properties of the stationary states] \label{theorem: poincare maps and fixed points: stability}
When $\xi>0.5$, the stationary states are: (i) asymptotically stable in the combinations of regions $(1,5)$, $(1,7)$, $(4,7)$, $(3,5)$ and $(3,8)$; (ii) Lyapunov stable in the combination of regions $(2,6)$; (iii) unstable in the combination of regions $(3,7)$. When $\xi< 0.5$, the stationary states are: (i) asymptotically stable in the combinations of regions $(1,5)$ and $(4,8)$; (ii) Lyapunov stable in the combination of regions $(2,6)$; (iii) unstable in the combinations of regions $(4,7)$ and $(3,8)$. When $\xi=0.5$, the stationary states are: (i) asymptotically stable in the combination of regions $(1,5)$; (ii) Lyapunov stable in the combinations of regions $(2,6)$, $(4,7)$ and $(3,8)$.
\end{theorem}
\textbf{Proof:} From Table \ref{table:possible poincare maps:perturbation}, the error term $\epsilon(t+T)$ can be written as $\epsilon(t+T)=\epsilon(t)e^c$, where $c$ is the coefficient. To obtain the stability properties, we only need to identify the values of the coefficient $c$. The fixed point is unstable if $c>0$, while it is stable if $c\leq 0$. Furthermore, the fixed point is asymptotically stable if $c<0$. Since $\xi_1=\xi_2=\xi$, we can have $\gamma_1=\gamma_4$, $\gamma_2=\gamma_5$. Normally, we assume $k_j\approx 5k_c$, so we have $\gamma_1>\gamma_5$ when $\xi>0.5$. In addition, we have $\gamma_3>\gamma_5$ when $\xi>0.5$, and $\gamma_3\leq \gamma_5$ when $\xi\leq 0.5$. Then the derivation of stability is easy and thus omitted here.\eop

\section{Analysis of macroscopic fundamental diagrams and gridlock times}\label{Macroscopic fundamental diagrams and gridlock time analysis}
\subsection{Macroscopic fundamental diagrams} \label{Macroscopic fundamental diagrams}
Once the stationary states under different retaining ratios are obtained, the MFD can be analytically derived. However, the calculation of the average network flow-rate in Equation \refe{equ:average flow-rate} involves the integral of out-fluxes $g_1(t)$ and $g_2(t)$, which is difficult to solve. Therefore, we use the following approximation when the signalized double-ring network has reached a stationary state, i.e., $k_1(nT)=k_1^{*}$, $n\in\mathbb{N}_0$.
\begin{eqnarray}
q(k)=&\frac{\int_{s=nT}^{(n+1)T}g_1(s)ds+\int_{s=nT}^{(n+1)T}g_2(s)ds}{2T}=\frac{2\int_{s=nT}^{(n+1)T}g_1(s)ds}{2T}\nonumber\\
\approx& \frac{\pi T (g_1(nT)+g_1(nT+\pi T))}{2T}=\pi\frac{g_1(k_1^{*})+g_1(k_1(nT+\pi T))}{2}.
\label{equ:average flow-rate:approx}
\end{eqnarray}

\begin{figure}[!h]
\centering
\subfigure[$0.5<\xi<1$]{\includegraphics[scale=0.2]{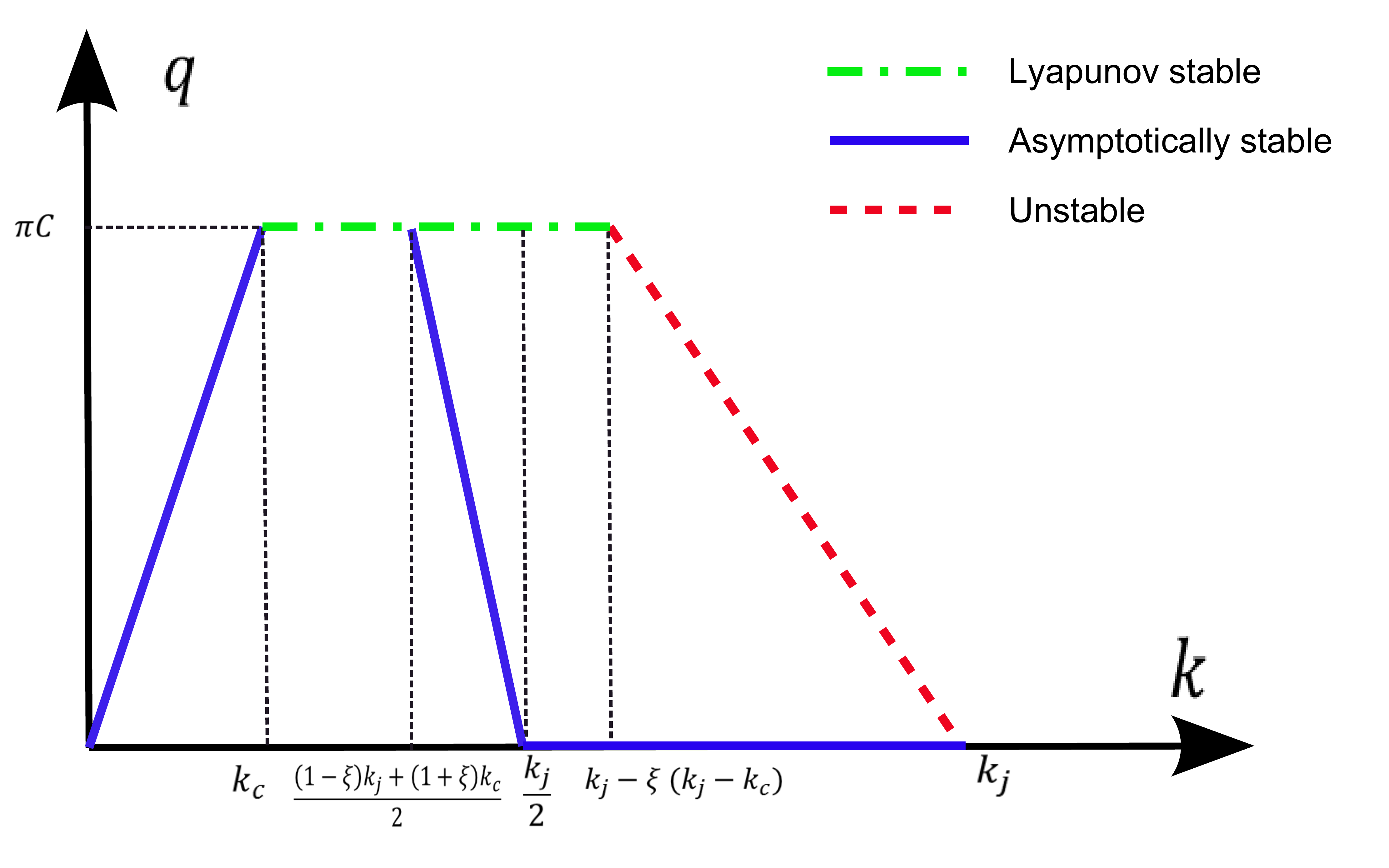}\label{fig:MFD-linkqueue:0.5greater}}
\subfigure[$ 0< \xi< 0.5$]{\includegraphics[scale=0.2]{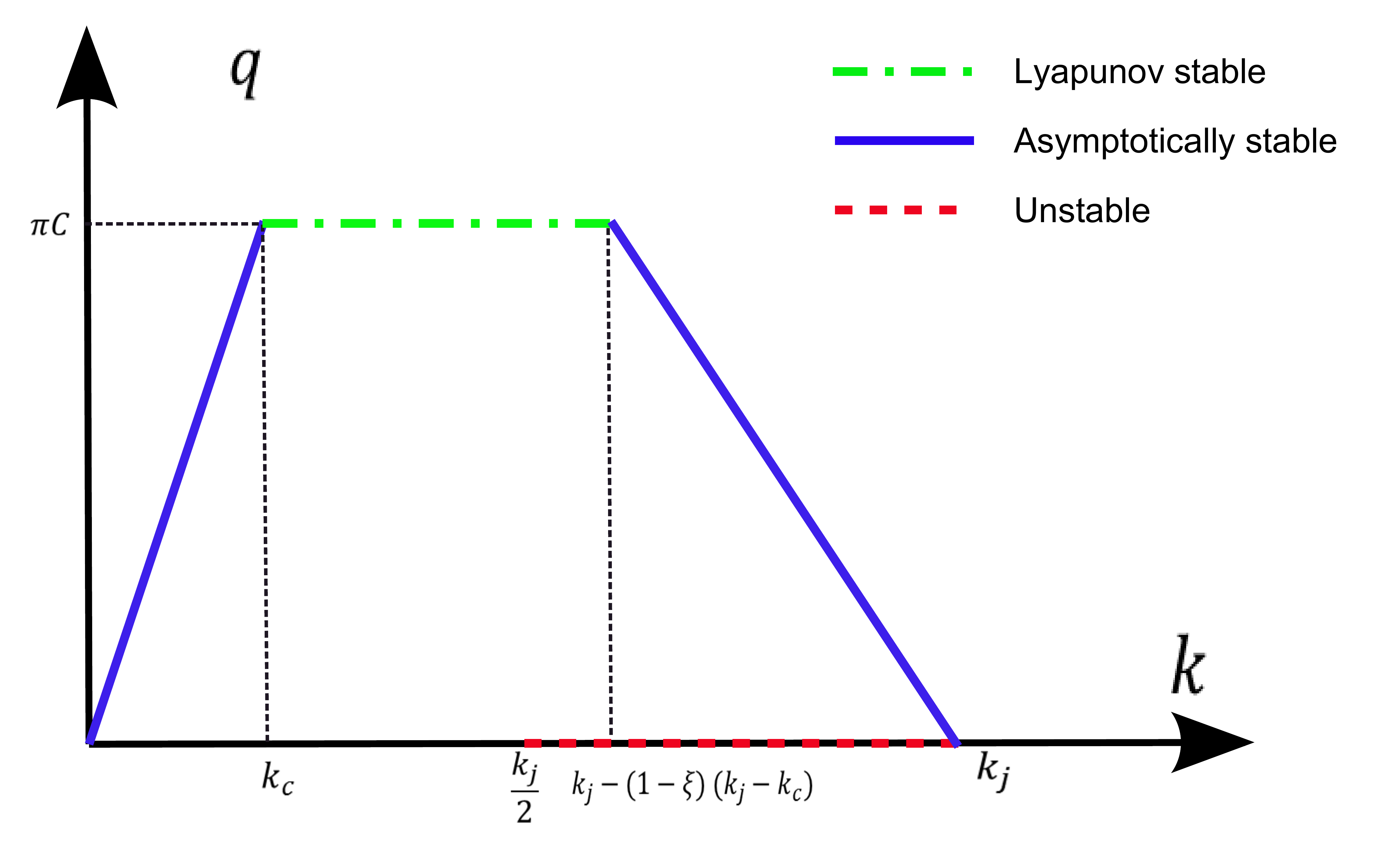}\label{fig:MFD-linkqueue:0.5smaller}}
\subfigure[$\xi=0.5$]{\includegraphics[scale=0.2]{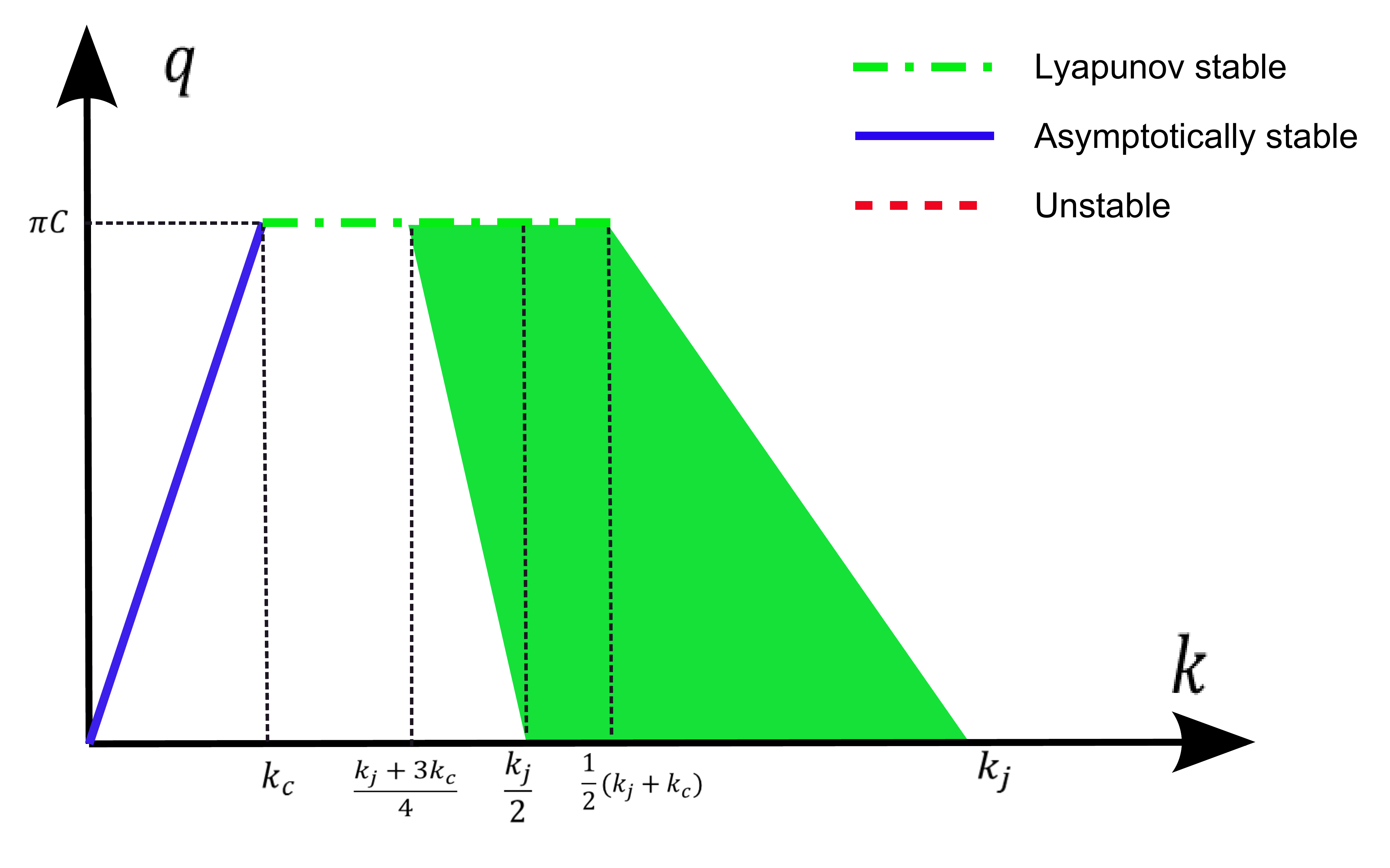}\label{fig:MFD-linkqueue:0.5}}
\caption{ MFDs for the signalized double-ring network with different retaining ratios. }
\label{fig:MFD-linkqueue}
\end{figure}

According to Equation \refe{equ:average flow-rate:approx}, we have the following MFD when $0.5<\xi<1$:
\begin{equation}
q(k;\pi,\xi)\approx \cas{{ll} \pi v_f k, & k\in[0, k_c]\\ \pi C, & k\in(k_c, k_j-\xi(k_j-k_c)] \\ \pi C\frac{k_j-2k}{\xi(k_j-k_c)-k_c}, & k\in(\frac{(1-\xi)k_j+(1+\xi)k_c}{2},k_j/2]\\\pi C\frac{k_j-k}{\xi(k_j-k_c)}, & k\in(k_j-\xi(k_j-k_c),k_j) \\ 0, & k\in(k_j/2,k_j]} 
\label{equ:MFD:greater:general}
\end{equation}
The proof is provided in \ref{appendix 1}. The MFD in Equation \refe{equ:MFD:greater:general} is shown in Figure \ref{fig:MFD-linkqueue:0.5greater}. From the figure, we find that: (i) multivaluedness exists in the network flow-density relation when $k\in(\frac{(1-\xi)k_j+(1+\xi)k_c}{2},k_j)$; (ii) for $k\geq \frac{k_j}{2}$, the network can have stationary states with non-zero flow-rates, but these stationary states are either Lyapunov stable or unstable; (iii) the gridlock states are asymptotically stable.

Similarly, when $0< \xi< 0.5$,  we have the following MFD:
\begin{equation}
q(k;\pi,\xi)\approx \cas{{ll} \pi v_f k, & k\in[0, k_c]\\ \pi C, & k\in(k_c, k_j-(1-\xi)(k_j-k_c)] \\\pi C\frac{k_j-k}{(1-\xi)(k_j-k_c)}, & k\in(k_j-(1-\xi)(k_j-k_c),k_j) \\ 0, & k\in(k_j/2,k_j]} 
\label{equ:MFD:smaller:general}
\end{equation}
The proof is similar to \ref{appendix 1} and thus omitted here. The MFD in Equation \refe{equ:MFD:smaller:general} is provided in Figure \ref{fig:MFD-linkqueue:0.5smaller}. From the figure, we find that: (i) the signalized double-ring network can maintain higher average network flow-rates when the retaining ratios are small; (ii) multivaluedness also exists when $k\geq k_j/2$; (iii) the gridlock states are unstable with small retaining ratios.

However, when $\xi=0.5$, we find that the initial densities in the combinations of regions $(2,6)$, $(4,7)$, and $(3,8)$ are all stationary states. Therefore, the MFD is very different from those in Equations \refe{equ:MFD:greater:general} and \refe{equ:MFD:smaller:general} and is provided below:
\begin{equation}
q(k;k_1,\pi,\xi)\approx \cas{{ll} \pi v_f k, & k\in[0, k_c]\\ \pi C, & k\in(k_c, \frac{k_j+3k_c}{4}] \\ \pi C\frac{k_j-2k+k_1}{(k_j-k_c)/2}, & k\in(\frac{k_j+3k_c}{4}, k_j] \text{ and } \\ & k_1\in[\max\{\frac{k_c(k_j-2k)}{k_j/2-3k_c/2}, 2k-k_j\}, \max\{2k-\frac{k_j+k_c}{2},k\}]}
\label{equ:MFD:0.5:general}
\end{equation}
The proof is similar to \ref{appendix 1} and thus omitted here. The MFD in Equation \refe{equ:MFD:0.5:general} is provided in Figure \ref{fig:MFD-linkqueue:0.5}. From the figure, we find that: (i) there are infinitely many stationary states when $k>\frac{k_j+3k_c}{4}$, and thus, the corresponding average network flow-rate covers the whole green region shown in the figure; (ii) the gridlock states are Lyapunov stable; (iii) there are no unstable stationary states.

\subsection{Gridlock times} \label{Gridlock time analysis}
In Theorem \ref{theorem: poincare maps and fixed points: stability}, we find that stationary states in the combinations of regions $(4,7)$ and $(3,8)$ are stable gridlock states when $\xi>0.5$. That means starting from any initial conditions in these regions, the signalized double-ring network will finally get gridlocked. The gridlock time can be analytically calculated from the Poincar\'{e} maps. For example, in the combination of regions $(3,8)$, we have the Poincar\'{e} map $k_1(t+T)=Pk_1(t)=k_j(1-e^{(\gamma_2-\gamma_3)\pi T})+k_1(t)e^{(\gamma_2-\gamma_3)\pi T}$. For an initial state $(k_1(0),k)$, we have 
\begin{equation}
k_1(T)=k_j(1-e^{(\gamma_2-\gamma_3)\pi T})+k_1(0)e^{(\gamma_2-\gamma_3)\pi T}.
\end{equation} 
Then after $n$ cycles, we have 
\begin{equation}
k_1(nT)=k_j-(k_j-k_1(0))e^{(\gamma_2-\gamma_3)n\pi T}.
\end{equation}
Since $\gamma_2<\gamma_3$, $k_1(nT)$ will converge to $k_j$ as $n\rightarrow +\infty$. If we define the gridlock time $T_g$ as the time when $k_1(T_g)\approx (1-\sigma) k_j$, where $\sigma$ is very small, we have
\begin{equation}
T_g\approx\frac{1}{\pi(\gamma_3-\gamma_2)} \ln\{\frac{k_j-k_1(0)}{\sigma k_j}\}.
\label{equ:time to gridlock}
\end{equation}

\begin{figure}[!h]
\centering
\subfigure[With different retaining ratios]{\includegraphics[scale=0.5]{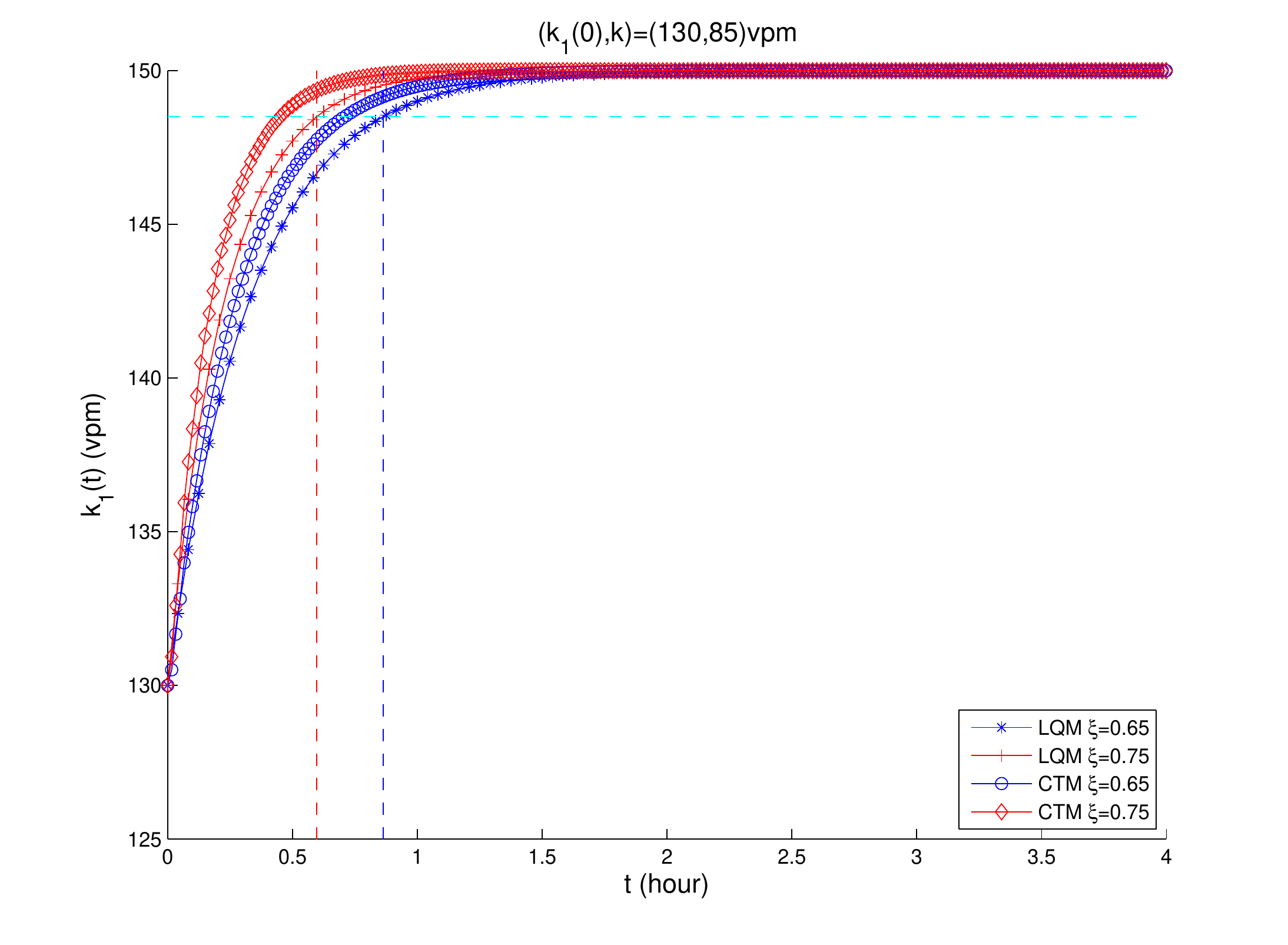}\label{fig:gridlock time:xi compared}}
\subfigure[With different initial densities]{\includegraphics[scale=0.5]{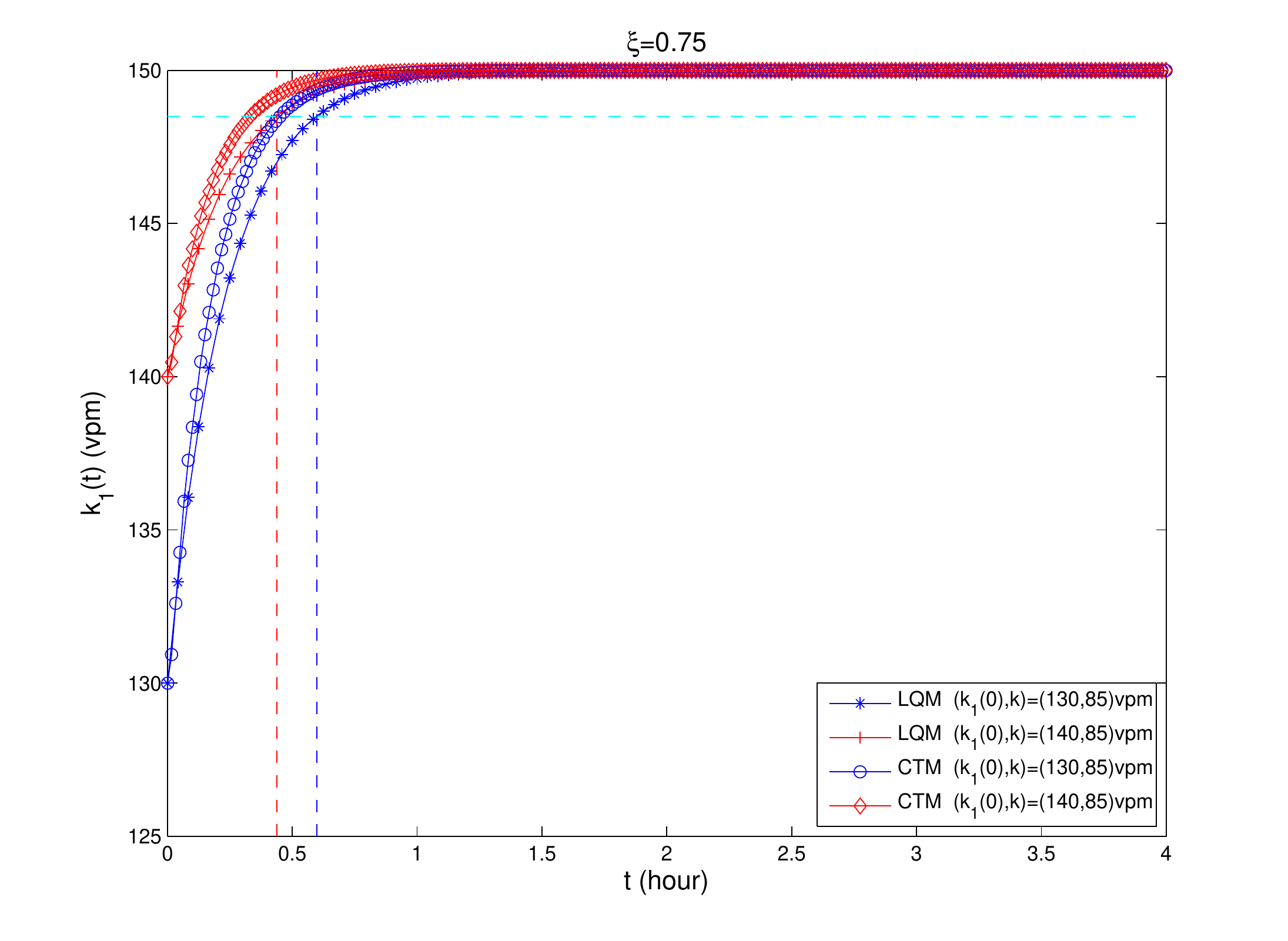}\label{fig:gridlock time:density compared}}
\caption{The gridlock patterns with different retaining ratios and initial densities.}
\label{fig:gridlock time}
\end{figure}

In Figure \ref{fig:gridlock time}, we provide the gridlock patterns with different retaining ratios and initial densities in the combination of regions $(3,8)$. In the figure, the horizontal dashed line is the threshold with $\sigma=0.01$ for gridlock conditions. The vertical dashed lines are the calculated gridlock times from Equation \refe{equ:time to gridlock}. The cycle length is 30s with a lost time of 2s for each phase. The updating time step $\Delta t$ is 0.01s in LQM simulations. The following two trends can be observed from the figure and verified by taking the derivatives with respect to $\xi$ and $k_1(0)$ in Equation  \refe{equ:time to gridlock}: (i) given the same initial densities, the network is harder to get gridlocked with lower retaining ratios, which is shown in Figure \ref{fig:gridlock time:xi compared}; (ii) given the same retaining ratios, the network will get gridlocked earlier if ring 1 is initially more congested, which is shown in Figure \ref{fig:gridlock time:density compared}. The gridlock patterns obtained from CTM simulations with $\Delta t=0.25$s are also provided in Figure \ref{fig:gridlock time}. We find that the two trends can also be found in CTM simulations, which indicates that they are the characteristics of the signalized double-ring network itself and are not related to the simulation models. Interestingly, we also find that with the same initial settings in Figure \ref{fig:gridlock time}, the network tends to get gridlocked earlier in the CTM.

\section{Numerical solutions to the Poincar\'{e} maps with long cycle lengths} \label{Numerical solutions to the Poincare maps with long cycle lengths}
With long cycle lengths, derivation of closed-form Poincar\'{e} maps is difficult since it is hard to track the time instants when the density evolution orbit crosses the boundaries of different regions. However, since the system dynamic equation (Equation \refe{equ:switched affine system}) and the coefficients of different regions (Table \ref{table:possible coefficients}) are known, it is possible to provide numerical methods to solve the Poincar\'{e} maps.  

Here, we denote $\Phi(k_1)=k_1-Pk_1$. For given average network density $k$, retaining ratio $\xi$, and signal settings $(\delta_1(t),\delta_2(t))$, finding the fixed points $k_1^{*}$ is the same as finding the roots in $\Phi(k_1)$. Note that the right-hand side of Equation \refe{equ:switched affine system} is a piecewise linear function and is continuous at the boundaries, and thus, $Pk_1$ as well as $\Phi(k_1)$ is continuous. But since the derivative of $\Phi(k_1)$ is not available, the secant method \citep{epperson2014introduction} is used and provided below: 
\begin{equation}
k_1^{n+1}=k_1^n-\Phi(k_1^n)[\frac{k_1^n-k_1^{n-1}}{\Phi(k_1^n)-\Phi(k_1^{n-1})}].
\label{equ:the secant method}
\end{equation}
$k_1^{n+1}$ is the updated density at step $n+1$. $k_1^{n}$, $k_1^{n-1}$, $\Phi(k_1^n)$ and $\Phi(k_1^{n-1})$ are the densities and function values at steps $n$ and $n-1$, respectively. 

\begin{figure}
\centering
\includegraphics[scale=0.8]{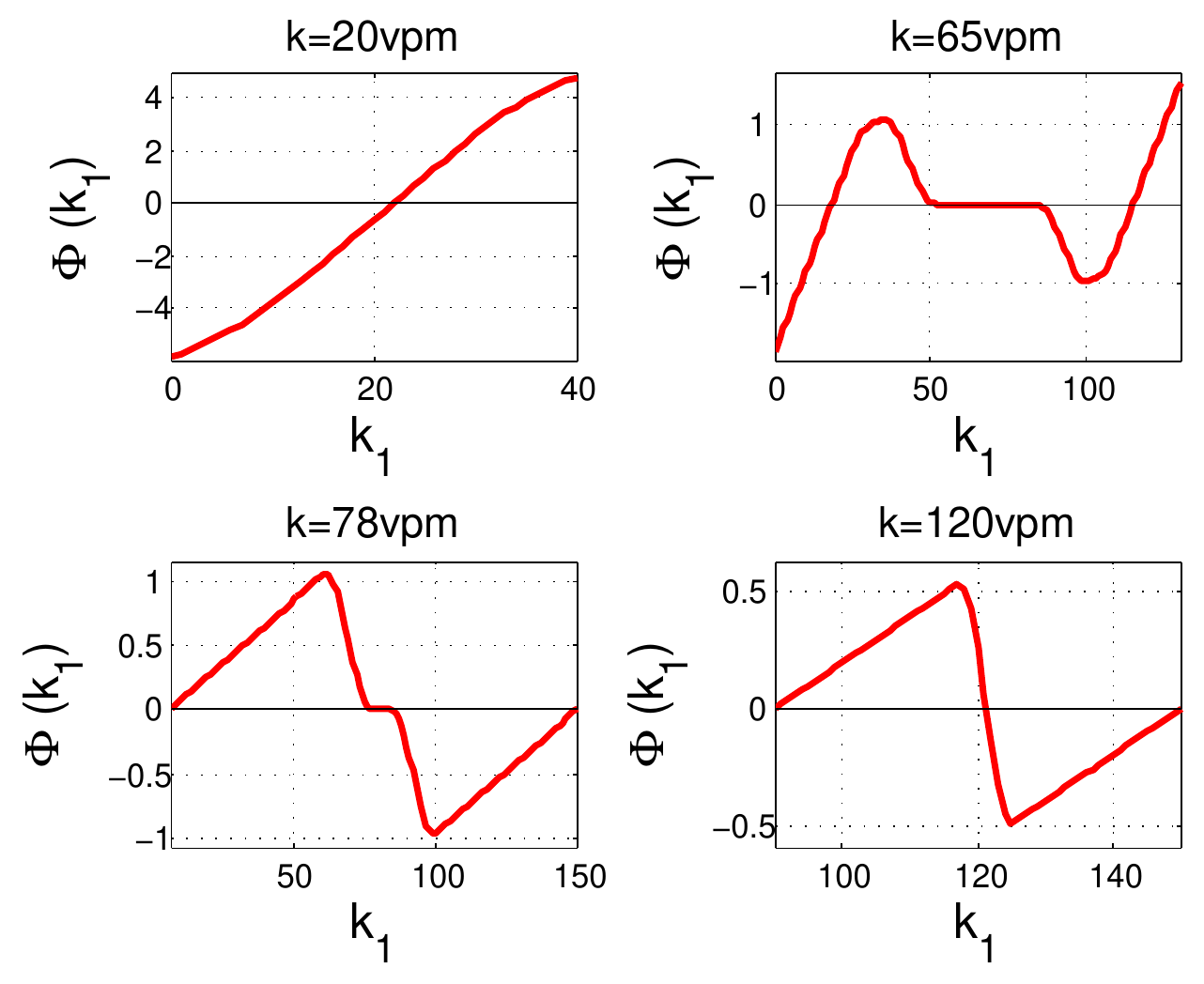}
\caption{Relations between $\Phi(k_1)$ and $k_1$ under different average network densities when $\xi=0.55$, $T=60$s, and $\Delta=4$s.}
\label{fig:relation between poincare map function value and ring 1 density}
\end{figure}

In Figure \ref{fig:relation between poincare map function value and ring 1 density}, we provide the relations between $\Phi(k_1)$ and $k_1$ under different average network densities when $\xi=0.55$, $T=60$s, and $\Delta=4$s. From the figure, we know that $\Phi(k_1)$ can have one root, multiple roots, or infinite roots under different average network densities. That means starting from different $k_1^{0}$, we may get different roots in the end using numerical methods. Therefore, a brute-force search in $k_1$ is needed to obtain a full map of the stationary states. In addition, when $\Phi(k_1)$ has infinite roots (e.g., when $k=65$vpm or $78$vpm), judgements on the values of $\Phi(k_1^{n})$ and $\Phi(k_1^{n-1})$ are needed to avoid zero values in the denominator in Equation \refe{equ:the secant method}. The algorithm of finding stationary states is provided in \ref{appendix 2}. 

\begin{figure}
\centering
\subfigure[$\xi=0.85$]{\includegraphics[scale=0.5]{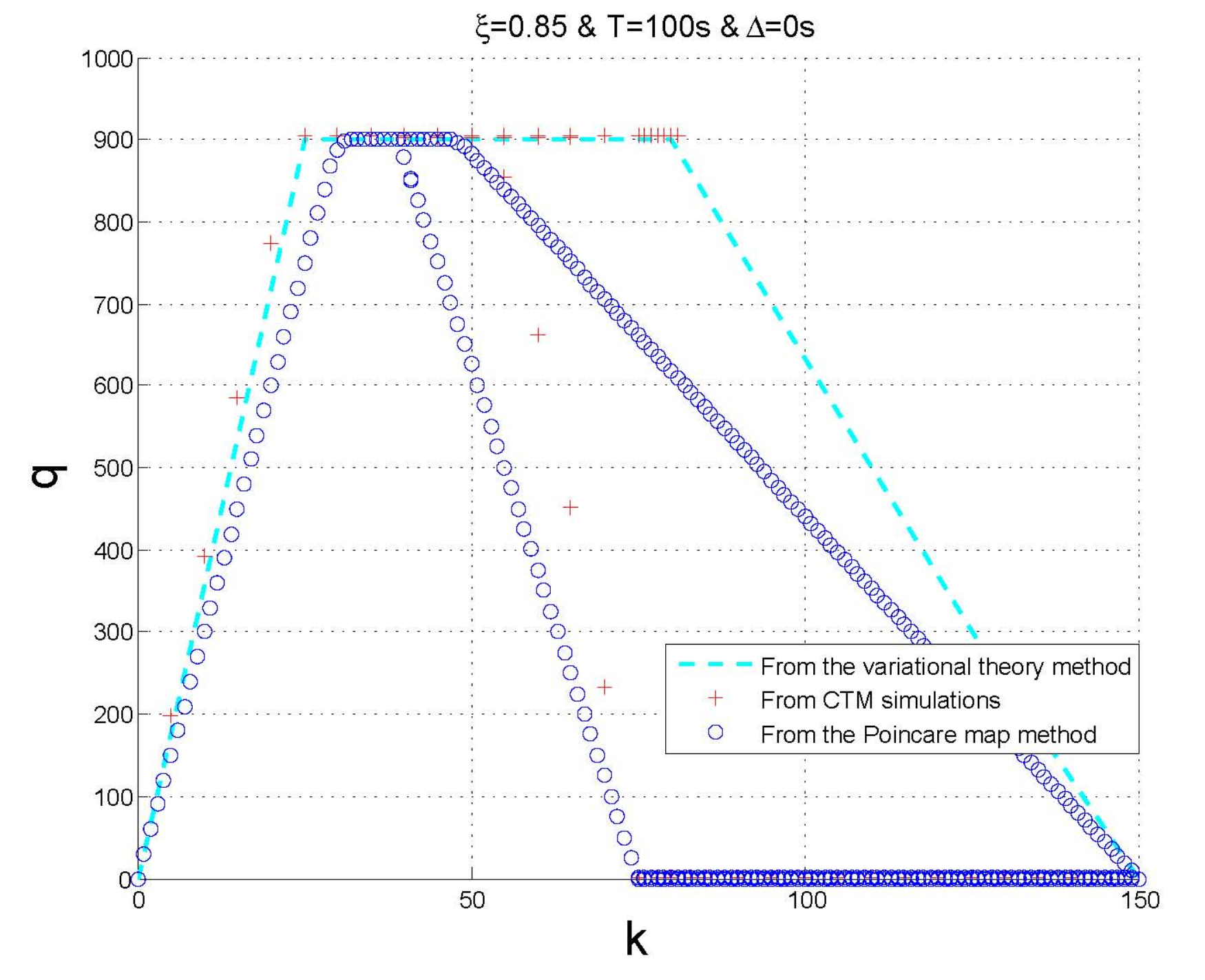}\label{fig:MFDs with diferent retaining ratios:85}}
\subfigure[$\xi=0.3$]{\includegraphics[scale=0.55]{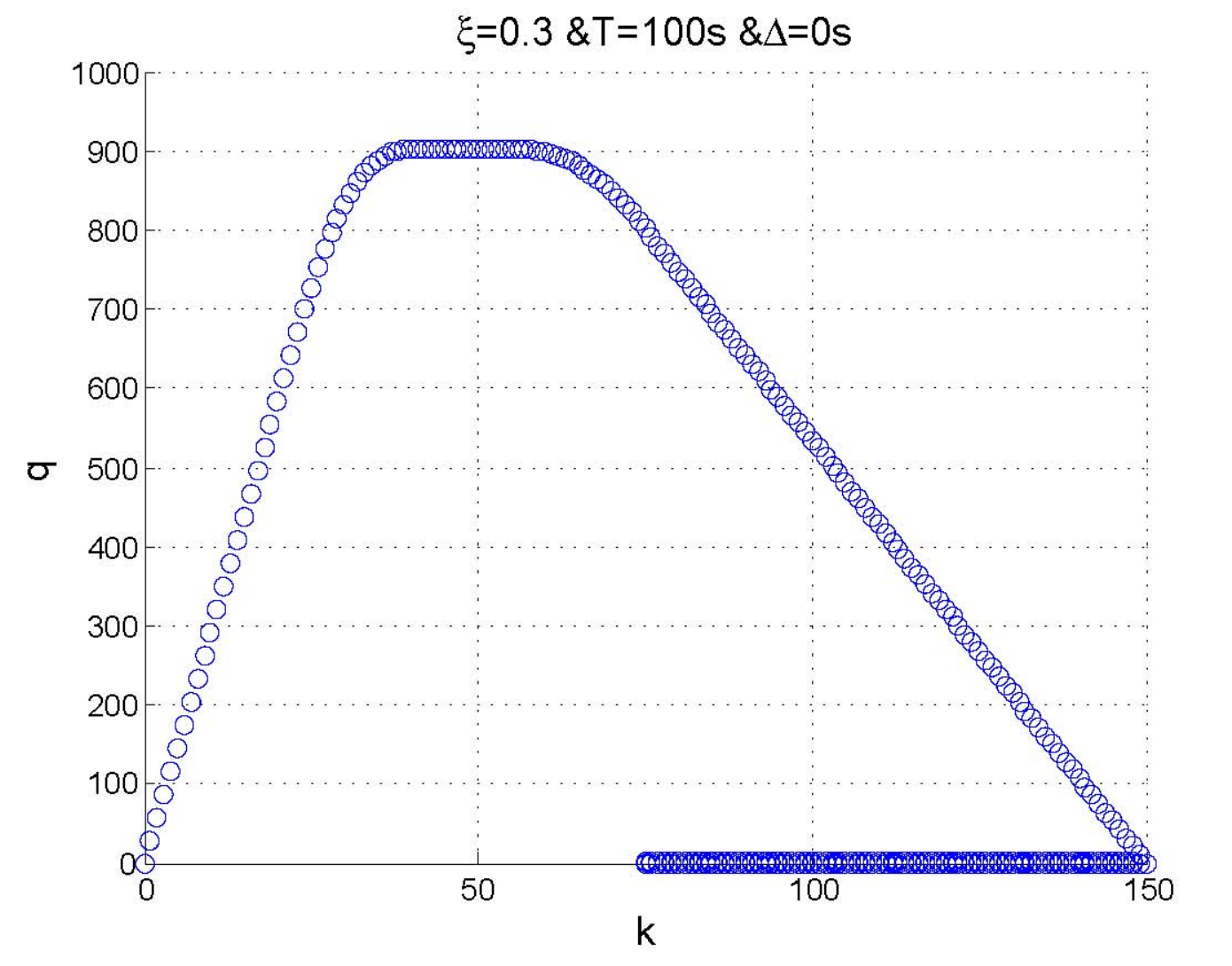}\label{fig:MFDs with diferent retaining ratios:30}}
\subfigure[$\xi=0.5$]{\includegraphics[scale=0.55]{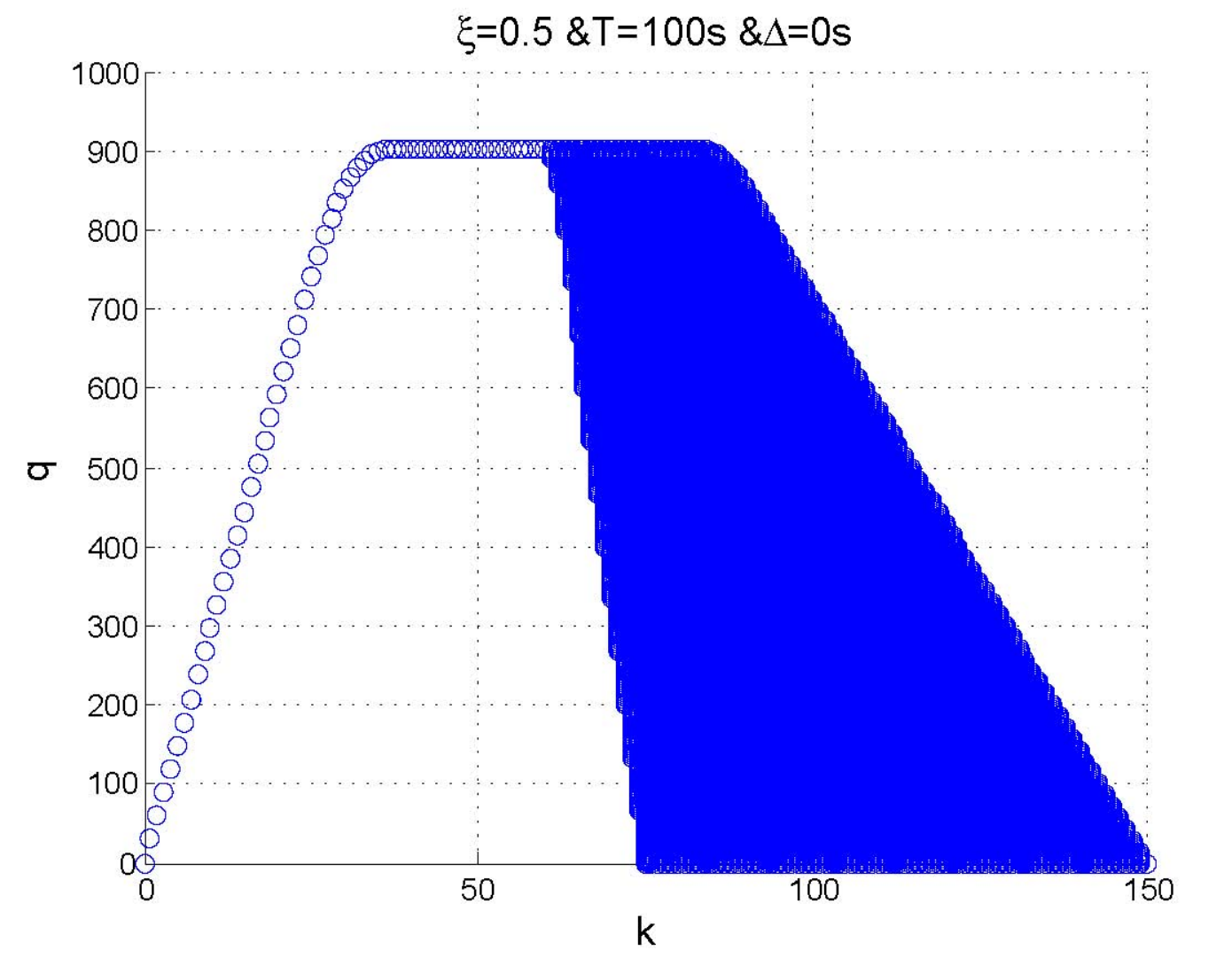}\label{fig:MFDs with diferent retaining ratios:50}}
\caption{MFDs with $T=100$s, $\Delta=0$s, and different retaining ratios.}
\label{fig:MFDs with diferent retaining ratios}
\end{figure}

In Figure \ref{fig:MFDs with diferent retaining ratios:85},  we provide the MFD in the signalized double-ring network with $\xi=0.85$, $T=100$s, and $\Delta=0$s using Equation \refe{equ:the secant method}. As a comparison, we also provide the MFDs obtained from the variational theory method and CTM simulations \citep{jin2013kinematic} with the same settings in the same figure. From the figure, we find that the MFDs obtained from the three methods are different. The variational theory method fails to observe the multivaluedness and gridlock phenomena in the MFD. Even though such two phenomena are observed using CTM simulations, the unstable congested branch with non-zero average network flow-rates is not observed. Using Equation \refe{equ:the secant method}, we can get a complete map of the MFD: the multivaluedness and gridlock phenomena as well as the unstable congested branch are all observed. The obtained MFD is similar to the one in Figure \ref{fig:MFD-linkqueue:0.5greater} with small cycle lengths. Since detailed traffic dynamics inside a link are aggregated at the link level in the link queue model, it is not surprised to find that the flow-rates obtained from the Poincar\'{e} map method are systematically lower than or equal to those obtained from the variational theory method and CTM simulations.

In addition, in Figures \ref{fig:MFDs with diferent retaining ratios:30} and \ref{fig:MFDs with diferent retaining ratios:50}, we provide the MFDs obtained from the Poincar\'{e} map method with $\xi=0.3$ and $\xi=0.5$, respectively. From the two figures, we find that the MFDs with long cycle lengths are very similar to the ones with short cycle lengths in Figures \ref{fig:MFD-linkqueue:0.5smaller} and \ref{fig:MFD-linkqueue:0.5}.

\section{Conclusions}\label{Conclusions}
In this paper, a signalized double-ring network was studied as a switched affine system under the framework of the link queue model \citep{jin2012link} and the assumption of a triangular traffic flow fundamental diagram \citep{haberman1977mathematical}. Stationary states with periodic density evolution orbits were shown to exist as the switched affine system periodically visits the $(k_1,k)$ space in Figure \ref{fig:phase_periodic_maps}. Poincar\'{e} maps were introduced to analyze the existence and the property of these stationary states. With short cycle lengths, closed-form Poincar\'{e} maps were obtained. Stationary states and their corresponding stability properties were derived and analyzed by finding and analyzing the fixed points on the Poincar\'{e} maps. It was found that under different combinations of signal settings, retaining ratios, and initial densities, stationary states can be Lyapunov stable, asymptotically stable, or unstable. Macroscopic fundamental diagrams were derived based on the stationary states, in which the network flow-rate is a function of the network density, route choice behaviors, and signal settings. In addition, the derived MFDs are more complete since the multivaluedness and gridlock phenomena as well as the unstable congested branch are all observed. Since the system will get gridlocked when $k\geq k_j/2$ and $\xi>0.5$, the gridlock time was also analyzed under different retaining ratios and initial densities. It was found that the network is harder to get gridlocked with lower retaining ratios, but it will get gridlocked earlier if one ring is initially more congested. Since closed-form Poincar\'{e} maps are hard to obtain with long cycle lengths, the secant method was used to numerically find the fixed points on the Poincar\'{e} maps. It was found that the obtained stationary states and the MFDs are very similar to those with short cycle lengths.

Different from the kinematic wave approach in \citep{jin2013kinematic} and the two bin model in \citep{daganzo2011macroscopic-b,gayah2011clockwise}, the introduction of the link queue model and the switched affine system enables us to analytically study the traffic statics and dynamics in the double-ring network with general signal settings. The way of using Poincar\'{e} maps to find the stationary states is physically meaningful and mathematically easier to solve. The findings in this paper, e.g., stationary states, stability, macroscopic fundamental diagrams, and gridlock patterns, are consistent with those in earlier studies \citep{daganzo2011macroscopic-b, gayah2011clockwise, jin2013kinematic}. But the solutions in this paper are more complete.

From our analysis, when $k\geq k_j/2$ and $\xi>0.5$, there exist two types of stationary states in the signalized double-ring network: the asymptotically stable gridlock states and the unstable states with more symmetric density distributions and higher average network flow-rates. To improve the network performance, we are interested in developing new signal control and route choice strategies to adaptively change the effective green times and retaining ratios. In addition, we are also interested in applying the analytical framework in this paper to solve the traffic dynamics in more general networks, e.g., a signalized grid network, in the future. 

\section*{Appendix A. Approximation of average network flow-rates}\label{appendix 1} 
The stationary states are provided in Table \ref{table:possible poincare maps} under different retaining ratios. For $\xi>0.5$, the possible combinations of regions having stationary states are: $(1,5)$, $(1,7)$, $(2,6)$, $(4,7)$, $(3,5)$, $(3,7)$, and $(3,8)$. For regions $(4,7)$ and $(3,8)$, the stationary states are gridlock states, and therefore, the average network average flow-rates are zero. For the rest of the regions, we can approximate the average network flow-rate using Equation \refe{equ:average flow-rate:approx}.
\begin{itemize}
\item [(1)] For regions $(1,5)$, the fixed point is $k_1^{*}=\frac{2k}{1+e^{-\gamma_1 \pi T}}$. Starting with $k_1(nT)=k_1^{*}$, we can get $k_1(nT+\pi T)=k_1(nT) e^{-\gamma_1 \pi T}$. Since ring 1 is uncongested, $g_1(k_1)=v_fk_1$. Therefore, the average network flow-rate is
\begin{eqnarray}
q(k)\approx \pi\frac{g_1(k_1^{*})+g_1(k_1(nT+\pi T))}{2}&=\frac{1}{2}\pi v_f (k_1^{*}+k_1(nT+\pi T))\nonumber\\
&=\frac{1}{2}\pi v_f 2k \frac{1+e^{-\gamma_1 \pi T}}{1+e^{-\gamma_1 \pi T}}=\pi v_f k.
\end{eqnarray}

\item [(2)] For regions $(1,7)$, the fixed point is $k_1^{*}=\frac{(k_j-2k)(e^{\gamma_5\pi T}-1)}{1-e^{(\gamma_5-\gamma_1)\pi T}}$. Starting with $k_1(nT)=k_1^{*}$, we can get $k_1(nT+\pi T)=k_1(nT) e^{-\gamma_1 \pi T}$. Since ring 1 is uncongested, $g_1(k_1)=v_fk_1$. In addition, since T is small, $-\gamma_1 \pi T$, $\gamma_5 \pi T$, and $(\gamma_5-\gamma_1) \pi T$ are also small. Therefore, the average network flow-rate is
\begin{eqnarray}
q(k)& \approx \frac{1}{2}\pi v_f (k_1^{*}+k_1(nT+\pi T))=\frac{1}{2}\pi v_f (k_j-2k) \frac{(e^{\gamma_5\pi T}-1)(1+e^{-\gamma_1 \pi T})}{1-e^{(\gamma_5-\gamma_1)\pi T}}\nonumber\\
&\approx\frac{1}{2}\pi v_f (k_j-2k) \frac{\gamma_5\pi T(2-\gamma_1 \pi T)}{(\gamma_1-\gamma_5)\pi T}\approx \pi C \frac{(k_j-2k)}{\xi(k_j-k_c)-k_c}.
\end{eqnarray}

\item [(3)] For regions $(2,6)$, the fixed point is $k_1^{*}=k_1(t)$. In this combination of regions, the out-fluxes are restricted by the capacity. Therefore, the average network flow-rate is 
\begin{equation}
q(k)\approx \pi\frac{g_1(k_1^{*})+g_1(k_1(nT+\pi T))}{2}=\pi C.
\end{equation}

\item[(4)] For regions $(3,5)$, the fixed point is $k_1^{*}=\frac{2k(1-e^{-\gamma_4\pi T})-k_j(e^{\gamma_2\pi T}-1)e^{-\gamma_4\pi T}}{1-e^{(\gamma_2-\gamma_4)\pi T}}$. Starting with $k_1(nT)=k_1^{*}$, we can get $k_1(nT+\pi T)=k_j(1-e^{\gamma_2\pi T})+k_1(nT)e^{\gamma_2\pi T}$. In this combination of regions, the out-flux is governed by the supply in ring 1, i.e., $g_1(k_1)=\frac{S_1(k_1)}{\xi}=\frac{C(k_j-k_1)}{\xi(k_j-k_c)}$. Therefore, the average network flow-rate is
\begin{eqnarray}
q(k)&\approx \pi C \frac{2k_j-(k_j(1-e^{\gamma_2\pi T})+k_1^{*}(1+e^{\gamma_2\pi T}))}{2\xi(k_j-k_c)}\approx \pi C \frac{2k_j-2\frac{2k*\gamma_4-k_j*\gamma_2}{\gamma_4-\gamma_2}}{2\xi(k_j-k_c)} \nonumber\\
&= \pi C \frac{(k_j-2k)\frac{\gamma_4}{\gamma_4-\gamma_2}}{\xi(k_j-k_c)} = \pi C \frac{(k_j-2k)}{\xi(k_j-k_c)-k_c}.
\end{eqnarray}

\item[(5)] For regions $(3,7)$, the fixed point is $k_1^{*}=\frac{2k+k_j(e^{\gamma_2\pi T}-1)}{e^{\gamma_2\pi T}+1}$. Starting with $k_1(nT)=k_1^{*}$, we can get $k_1(nT+\pi T)=k_j(1-e^{\gamma_2\pi T})+k_1(nT)e^{\gamma_2\pi T}$. In this combination of regions, the out-flux is governed by the supply in ring 1, i.e., $g_1(k_1)=\frac{S_1(k_1)}{\xi}=\frac{C(k_j-k_1)}{\xi(k_j-k_c)}$. Therefore, the average network flow-rate is
\begin{eqnarray}
q(k)&\approx \pi C \frac{2k_j-(k_j(1-e^{\gamma_2\pi T})+k_1^{*}(1+e^{\gamma_2\pi T}))}{2\xi(k_j-k_c)}= \pi C \frac{2k_j-(k_j(1-e^{\gamma_2\pi T})+2k+k_j(e^{\gamma_2\pi T}-1))}{2\xi(k_j-k_c)} \nonumber\\
&=\pi C \frac{k_j-k}{\xi(k_j-k_c)}.
\end{eqnarray}
\end{itemize}

Based on the average network flow-rates calculated above, the macroscopic fundamental diagram for the signalized double-ring network with $0.5<\xi<1$ can be easily derived and thus omitted here.
\eop

\section*{Appendix B. Algorithm of finding stationary states} \label{appendix 2}
\begin{itemize}
\item[] \textbf{Inputs:} $k$, $T$, $\Delta$, $\xi$, and $\Phi(k_1)$
\item[] \textbf{Initialization:} vector of stationary states, i.e., SS=[]; minimum value of $k_1$, i.e., $k_{1,min}=\max\{2k-k_j,0\}$; maximum value of $k_1$, i.e., $k_{1,max}=\min\{2k,k_j\}$; the threshold of $k_1$, i.e., $e_k$; searching step, i,e., $\Delta k$; maximum number of iterations, i.e., $n_{max}$
\item[] \textbf{For} $k_1=k_{1,min}:\Delta k:k_{1,max}$
\begin{itemize}
\item[] Set $k_1^{0}=k_1$ and calculate $\Phi(k_1^{0})$ 
\item[] \textbf{If} $\Phi(k_1^{0})==0$
\begin{itemize}
\item[] $k_1^{0}$ is a root, and add it to SS
\end{itemize}
\item[] \textbf{Else}
\begin{itemize}
\item[] Set $k_1^{1}=Pk_1^{0}$ and calculate $\Phi(k_1^{1})$ 

\item[] \textbf{For} n=$1:n_{max}$
\begin{itemize}
\item[] \textbf{If} $|k_1^{1}-k_1^{0}|<e_k$ \textbf{or} $\Phi(k_1^{1})==0$, add $k_1^{1}$ to SS and break
\item[] $k_1^{tmp}=k_1^1-\Phi(k_1^1)[\frac{k_1^1-k_1^{0}}{\Phi(k_1^1)-\Phi(k_1^{0})}]$, and calculate $\Phi(k_1^{tmp})$
\item[] $k_1^0=k_1^1$ and $\Phi(k_1^{0})=\Phi(k_1^{1})$
\item[] $k_1^1=k_1^{tmp}$ and $\Phi(k_1^1)=\Phi(k_1^{tmp})$
\end{itemize}
\item[] \textbf{End for}
\end{itemize}
\item[] \textbf{End if}
\end{itemize}
\item[] \textbf{End for}
\end{itemize}

\end{document}